\documentclass[a4paper,11pt,bibliography=totoc]{scrartcl}

\usepackage[USenglish]{babel}
\usepackage[T1]{fontenc}
\usepackage[utf8]{inputenc}

\usepackage{graphicx,subfigure,tikz,pgfplots}
\usepackage[customcolors]{hf-tikz}
\hfsetbordercolor{white}
\definecolor{myblue}{RGB}{17,64,111}
\definecolor{myred}{RGB}{150,30,15}
\pgfplotsset{compat=1.17}

\usepackage{amsmath,amssymb,amsfonts,amsthm,array,paralist,todonotes,mathtools}
\usepackage{stmaryrd}
\SetSymbolFont{stmry}{bold}{U}{stmry}{m}{n}
\usepackage{array,tablefootnote,thmtools}
\usepackage[toc,page]{appendix}

\usepackage[pdftex,colorlinks=true,linkcolor=blue,citecolor=green,urlcolor=blue,bookmarks]{hyperref}

\newcommand{\R}{\mathbb{R}}

\newcommand{\dd}{\textup{d}}
\newcommand{\vol}{\text{vol}}

\newcommand{\bb}[1]{\mathbf{#1}}
\newcommand{\bs}[1]{\boldsymbol{#1}}
\newcommand{\dbl}{\left\llbracket}
\newcommand{\dbr}{\right\rrbracket}
\def\c{{\textup c}}
\def\p{{\textup p}}
\def\s{{\textup s}}
\def\t{{\textup t}}
\def\x{{\mathbf x}}
\def\y{{\mathbf y}}
\def\z{{\mathbf z}}
\def\xin{{\textup x}}
\def\yin{{\textup y}}
\def\zin{{\textup z}}

\newtheorem{thm}{Theorem}[section]

\newenvironment{prf}{\textbf{Proof:}} {\hspace*{\fill} $\square$ \newline}

\begin{document}
     
    \title{Balance Laws and Transport Theorems for Flows with Singular Interfaces}
    \author{Ferdinand Thein\footnotemark[1],\; Gerald Warnecke\footnotemark[2]}
    \date{\today}
    \maketitle
		\noindent \small
Dedication: This paper is dedicated to the memory of Wolfgang Dreyer. For both authors he was a very important
source of knowledge as well as a partner for fruitful discussions on the topics of continuum mechanics and thermodynamics,
especially concerning phase transitions.\\[10mm]

\normalsize
\textbf{Abstract:} This paper gives a concise but rigorous mathematical description of a material control volume that is
separated into two parts by a singular surface at which physical states are discontinuous. 
The geometrical background material is summarized in a unified manner.
Transport theorems for use in generic balance laws are given with proofs since they provide some insight into the results.
Also the step from integral balances to differential equations is given in some detail. 
    \renewcommand{\thefootnote}{\fnsymbol{footnote}}
    \footnotetext[1]{Institut f\"ur Mathematik, Johannes-Gutenberg-Universit\"at Mainz,
     Staudingerweg 9, 55128 Mainz, Germany.
    \href{mailto:fthein@uni-mainz.de}{{\em fthein@uni-mainz.de}}}
    \footnotetext[2]{Institute for Analysis and Numerics, Otto-von-Guericke University Magdeburg, PSF 4120, D-39016 Magdeburg, Germany.
    \href{mailto:gerald.warnecke@ovgu.de}{{\em gerald.warnecke@ovgu.de}}}
    \renewcommand{\thefootnote}{\arabic{footnote}}



  %
	%
  \section{Introduction}

	In this work we present a general mathematical derivation of balance laws involving curved interfaces. These may be shock fronts, 
	phase boundaries, reaction fronts or surfaces containing surfactants. This includes
	the possibility of modeling physical processes, such as phase transitions or chemical reactions, and flow of
	substances within the interface. We want to provide a compact, self-contained reference for mathematical 
	modeling with hyperbolic balance laws, mathematical analysis and numerical computations of engineering 
	applications using such balance laws. We consider transport theorems, jump conditions as well as the 
	influence of curvature of the surfaces. 
	
	The main aim of this paper is a rigorous mathematical derivation of
	some general surface differential equations in Subsection \ref{subsec:jump}
	that combine differences in bulk quantities with dynamics on the surface.
	We give a unified treatment of volume and surface balances from a mathematical point of view 
	that draws on multiple sources in the literature. They
	each only provide partial aspects of our exposition. Also these sources are not easily compared due to differences in 
	assumptions, notation and presentation. We will show how they are related.
	An important point is the clarification of various derivatives, that are in use, of quantities
	defined on surfaces. Our proofs make heavy use of parametrizations of the surfaces.
	But, the results that we obtain are independent of the choice of parametrization.
	
	Our work was motivated by a preprint
	of Wolfgang Dreyer \cite{Dreyer2003} that has been included in this book. This again was based on M\"uller \cite{Mueller1985}.
	We focus on a detailed mathematical description of moving surfaces and the derivation of generic balance laws 
	in the volume as well as on the surface. This includes some proofs of fundamental theorems which we include because
	they give insight on how certain terms in the resulting differential equations are generated. For example in the derivation of
	M\"ueller \cite[(3.7)]{Mueller1985} a flux correction is postulated. In the proof of the generalized Reynolds
	Transport Theorem \ref{thm:gen_reynolds_transport} one can see how this correction arises from the description
	of a moving surface by parametrization and flow maps. 	
	Further sources of background material are Aris \cite{Aris1989}, Cermelli et al.\ \cite{Cermelli2005}, 
	Dziuk and Elliott \cite{Dziuk2007}, Grinfeld \cite{Grinfeld2013}, Gurtin \cite{Gurtin1993}, 
	Petryk and Mroz \cite{Petryk1986} as well as Slattery et al.\ \cite{Slattery2007}.
		
	The dynamics of the fluid motion can be described by certain physical laws. These laws are in general referred 
	to as {\em balance laws}. They are formulated by applying basic physical principles to moving material
	volumes and surfaces of arbitrary shape and size. Transport theorems and mathematical arguments are needed to derive
	differential equations from the integral balances. 
	Our derivation of differential equations on volumes and surfaces
	is based on Dreyer \cite{Dreyer2003} and M\"uller \cite{Mueller1985}.
	We also recommend Aris \cite{Aris1989} as well as Truesdell and Toupin \cite{Truesdell1960}.
	The latter is also very interesting from a historical perspective. 
	Truesdell and Toupin \cite{Truesdell1960} present the results in the context of original works.
	Thus the influences of mathematicians like {\em Christoffel, Euler, Hadamard, Hugoniot, Reynolds} on this topic become visible.
	A treatment of this topic in two dimensions can be found in Gurtin \cite{Gurtin1993}.
	Truesdell and Toupin \cite{Truesdell1960} did not treat dynamics on surfaces. Dziuk and Elliott \cite{Dziuk2007}
	only considered moving hypersurfaces in $\R^{n+1}$ and no volume dynamics. Another very useful reference 
	on interfacial transport is Slattery et al.\ \cite{Slattery2007}.
	
	We consider here the transport of control volumes and moving surfaces that remain inside some domain $\Omega\subseteq\R^3$.
	An important aspect that we do not take into account is the additional complication that arises when these 
	interact with the boundary $\partial\Omega$ of a bounded domain. The interaction of the surface with
	the boundary of a domain was studied 
	by Dell'Isola and Romano \cite{Dellisola1986}, Gurtin et al.\ \cite{Gurtin1989} as well as Estrada and Kanwal \cite{Estrada1991}.
		
	From a mathematical point of view balance laws are often conservation laws in the sense that the differential equations are
	given in divergence form.
	The field of conservation laws, even when restricted to particular equations, is far too extensive to be completely discussed here.
	For a brief overview of the mathematical results for balance laws, in particular {\em hyperbolic conservation laws}, 
	we refer to the following literature. Evans \cite{Evans1998} contains a compact description of the topic, 
	especially of systems of conservation laws.
	Further we recommend Dafermos \cite{Dafermos2016}, LeFloch \cite{LeFloch2002} and Smoller \cite{Smoller2012}.
	Especially Dafermos \cite{Dafermos2016} starts with a very general approach, introducing the balance laws in terms of suited 
	measures and further gives a nice historical review of the topic.
	Smoller \cite{Smoller2012} does not focus on hyperbolic conservation laws alone, reaction-diffusion equations are also treated.
	An introduction and discussion of certain analytical methods for conservation laws can be found in Warnecke \cite{Warnecke1999}.
	Some more recent analytical results were presented in Feireisl et al.\ \cite{Feireisl2021}.
	For additional literature see the references in the literature just mentioned.
	
	In this paper we do not provide specific balances for physical quantities. The latter can be found in 
	Bedeaux \cite{Bedeaux1986}, Dreyer \cite{Dreyer2003}
	M\"uller \cite{Mueller1985} as well as Truesdell and Toupin \cite{Truesdell1960}.	
	Mavrovouniotis and Brenner \cite{Mavrovouniotis1993} discussed a micro-macro approach to deriving
	surface balances as surface excess properties, a concept introduced by Gibbs. Further derivations of balance laws in the presence of 
	interfaces may be found in Anderson et al.\ \cite{Anderson2007}, Bedeaux et al.\
	\cite{Bedeaux1976}, Dell'Isola and Romano \cite{Dellisola1987} as well as Slattery et al.\ \cite{Slattery2007}.
		
	The key contribution from physics to this topic is provided by the theory of irreversible processes and 
	non-equilibrium thermodynamics. Some theory of non-equilibrium thermodynamics is given in 
	Bedeaux \cite{Bedeaux1986}, Mauri \cite{Mauri2013} as well as DeGroot and Mazur \cite{DeGroot2013}.
	An introduction to balance laws from a thermodynamic point of view can be found in Landau and Lifschitz\ \cite{Landau1987e}, 
	M\"uller \cite{Mueller1985} as well as M\"uller and M\"uller \cite{Mueller2009}.
	
	We begin by introducing in Section \ref{sec:calc_surf} various concepts for stationary surfaces. 
	This is then followed by the extension to moving surfaces. In the next section we formulate generic 
	integral balance laws that can be used for physical states such as mass, momentum or energy.
	Section \ref{sec:transport} contains various integral theorems for the volume and for surfaces that
	are needed in order to derive differential equations from the integral balances.
	These include balances for singular moving surfaces in the volume domain at which physical states may be discontinuous.
	A final section contains some comments on the lower dimensional cases.

	\section{Calculus of Stationary and Moving Surfaces}
	\label{sec:calc_surf}

	The aim of the present section is to briefly summarize the most important 
	results for moving surfaces. We refer to the works by Aris \cite{Aris1989}, Dziuk and Elliott \cite{Dziuk2007}, 
	Grinfeld \cite{Grinfeld2013} and the references therein. But, first we have to introduce various
	concepts for stationary surfaces.
	
	We will make extensive use of tensors. For a detailed introduction refer to the literature, 
	e.g.\ Aris \cite[Chapter 7]{Aris1989} or Grinfeld \cite[Chapter 6]{Grinfeld2013}. A tensor is a 
	quantity which has certain properties under the tranformation of the coordinate system. The order 
	of a tensor is given by the number of indices. A tensor with upper indices is called {\em contravariant} 
	and with lower indices {\em covariant}. A tensor may also be of mixed type with upper and lower indices.
	The indices in a term consisting of indexed quantities that appear only once are called {\em free indices}. They should
	appear in each term on each side of an equation exactly once.
	
	In the following we will use the following notational conventions:
	\begin{enumerate}[(i)]
	    \item Lower case Greek indices refer to the surface and may take the values $\{1,2\}$. 
			In contrast lower case Latin indices refer to the ambient space, i.e.\ the standard $\R^3$ in our case.
	    \item Bold lower case letters refer to vectors or tensors of order $1$ in $\R^3$, bold upper case letters to matrices
			or tensors of order $2$.
	    \item We use the Einstein summation convention, i.e.\ summation applies to every index which appears twice. 
			Once as superscript and once as subscript.
	    Such an index is often called {\em dummy index} and may be renamed. See the following to example of the scalar product
			in two and three space dimensions
	    \begin{align*}
					&\bb a\cdot\bb b =a^1b_1 + a^2b_2 = a^\alpha b_\alpha = a^\beta b_\beta,\\
	        &\mathbf{x}\cdot\mathbf{y} = \sum_{i=1}^3 \xin_i\yin_i = \xin^i\yin_i = \xin^j\yin_j
	    \end{align*}
			We will later also apply the convention to some summations from $1$ to $n$ using lower case Latin letters.
	    We call the summation of a dummy index a {\em contraction}.
	    Tensors may change in a certain way when transforming the coordinate system. A contraction is  
			invariant under the transformation and thus it is one of the fundamental operations in tensor calculus,
			see Grinfeld \cite[Chapter 6]{Grinfeld2013}. Also the contraction of tensors produces again a tensor.
	\end{enumerate}
	\subsection{Stationary Surfaces}
	\label{subsec:surf}

	A stationary two dimensional {\em surface} $\Sigma$ embedded in the ambient space $\R^3$ may be described through a 
	parametrization depending on a two dimensional coordinate system. We want to specify this in the following.
	Let $\mathcal{U} \subseteq \R^2$ as well as $\Omega\subseteq\R^3$ be open subsets, 
	$\bb u\in\mathcal{U}$ and $\bs\Phi : \mathcal{U} \to \Omega$ a bijective smooth
	parametrization such that
	\begin{align}
	    \Sigma \equiv \bs\Phi(\mathcal{U}) = \left\{\mathbf{x} \in \Omega\, |\, \mathbf{x} 
			= \bs\Phi (\bb u)=(\Phi^1(\bb u),\Phi^2(\bb u),\Phi^3(\bb u))^T \;\text{for}\;\bb u\in \mathcal{U}\right\}.
			\label{def:surface}
	\end{align}
	By a {\em smooth} function on an open set we will always mean that it is continuously differentiable on the set as 
	often as is necessary for formulas to be meaningful. It may also imply continuous bounded extensions to the boundary of the set
	when values on the boundary are needed. 
	
	Alternatively, such a surface can be defined implicitly as the level set of a
	smooth function $f:\Omega\to\R$. We will introduce such a description more generally 
	for moving surfaces in Subsection \ref{subsec:move}. 
	This is useful in integration theory, see e.g.\ Thomas et al.\ \cite[Section 16.5]{Thomas2004}, and
	in some applications, especially when dealing with a changing topology of 
	moving surfaces. A drawback is that for such a description of a surface one can only
	determine a normal velocity vector. In contrast, an explicit description of a moving surface determines via \eqref{eq:surf_velocity}
	a velocity vector of the surface. This explicit description will be a time dependent
	extension of the one given above.
	
	On the surface we must consider three kinds of quantities, such as tensors. There are those that are defined on the surface 
	only and are described by the surface coordinates $\bb u$ and depend only on them. Then there are those, 
	such as tangent and normal vectors, that have coordinates in the ambient space but depend on the surface coordinates.
	Further, there are those that have extensions in a neighborhood around the surface. Their nature determines what kind of
	operations are possible with or on them.
	
	The {\em tangent vector fields} $\bs\tau_\alpha\in\R^3$ and the {\em unit normal vector field} $\bs\nu\in\R^3$ 
	to the surface $\Sigma$ embedded in $\R^3$ are given by
	\begin{align}
	    \bs{\tau}_\alpha = \frac{\partial\bs\Phi}{\partial u_\alpha},\,\alpha = 1,2
	    \quad\text{and}\quad
	    \bs{\nu} = \frac{\bs{\tau}_1\times\bs{\tau}_2}{|\bs{\tau}_1\times\bs{\tau}_2|}.
			\label{def:tangent_normal}
	\end{align}
	Since we want the tangents to be linearly independent the Jacobian of $\bs\Phi$ should be of rank two.
	The components of the {\em surface metric tensor}\footnote{Aris \cite[Sections 9.13]{Aris1989} denoted this tensor by $a_{\alpha\beta}$
	and introduced an additional the metric tensor $g_{ij}$ of $\R^3$ in an exercise. For simplicity we assume the
	usual Euclidean scalar product on $\R^3$ where this tensor corresponds to the identity matrix. In \cite[Chapter 10]{Aris1989}
	he used such an extra metric tensor in an intrinsic formulation of the surface in $\R^2$.} are given as
	\begin{align}
	    g_{\alpha\beta} = \tau_{\alpha;i}\tau_{\beta;i}=\bs{\tau}_\alpha\cdot\bs{\tau}_\beta,
			\label{def:surface_metric}
	\end{align}
	cp.\ Aris \cite[(9.13.3)]{Aris1989} or Grinfeld \cite[(5.7)]{Grinfeld2013}.
	Thus the metric tensor $g_{\alpha\beta}$ is positive definite and symmetric.
	The inverse of the metric tensor is denoted by $g^{\alpha\beta}$. With the {\em Kronecker delta tensor} satisfying
	$\delta^\alpha_\gamma = 1$ for $\alpha =\gamma$
	and $=0$ for $\alpha \ne\gamma$ we have
	\begin{align}
	    g^{\alpha\beta}g_{\beta\gamma} = \delta^\alpha_\gamma\quad\text{and}\quad g 
			= \det(g_{\alpha\beta}) = g_{11}g_{22} - g_{12}^2.
			\label{def:inverse_det}
	\end{align}

	We consider the $\epsilon_{ijk}$ or {\em permutation tensor} for the cross product and determinant, 
	see Aris \cite[(2.32.4)]{Aris1989} or Grinfeld \cite[(9.8)]{Grinfeld2013}. 
	Its components are $\epsilon_{123}=\epsilon_{231}=\epsilon_{312}=1$,
	$\epsilon_{132}=\epsilon_{321}=\epsilon_{213}=-1$ and all other components equal $0$. So we have $1$ 
	for the even cyclic permutations of $123$,
	$-1$ for the odd permutations and $0$ for all repeated indices like $\epsilon_{112}$. This gives 
	$(\bb a \times\bb b)_i	=\epsilon_{ijk}a_jb_k$ and $\det\bb A =\epsilon_{ijk}a_{jk}$ for a matrix
	$\bb A =\left(a_{ij}\right)_{1\le i,j\le3}$. 	With the $\epsilon$-$\delta$-identity
	\begin{equation}
	\label{eq:ep_delta}
	\epsilon_{ijk}\epsilon_{irs}=\delta_{jr}\delta_{ks}-\delta_{js}\delta_{kr},
	\end{equation}
	see e.g.\ Aris \cite[Exercise 2.32.1]{Aris1989}, we have
	\begin{align*}
	|\bs{\tau}_1\times\bs{\tau}_2|^2=&\epsilon_{ijk}\tau_{\alpha;j}\tau_{\beta;k}\epsilon_{irs}\tau_{\alpha;r}\tau_{\beta;s}
	=(\delta_{jr}\delta_{ks}-\delta_{js}\delta_{kr})\tau_{\alpha;j}\tau_{\beta;k}\tau_{\alpha;r}\tau_{\beta;s}\\
	=&\tau_{\alpha;j}\tau_{\beta;k}\tau_{\alpha;j}\tau_{\beta;k}-\tau_{\alpha;j}\tau_{\beta;k}\tau_{\alpha;k}\tau_{\beta;j}
	=g_{\alpha\alpha}g_{\beta\beta}-(g_{\alpha\beta})^2 =g.
	\end{align*}
	This gives
	\begin{align}
	\label{eq:surface_elem}
	    |\bs{\tau}_1\times\bs{\tau}_2| =\sqrt{g}.
	\end{align}
	Let $\psi:\Sigma\to\R$ be a suitable function such that $\psi\circ\bs\Phi$ is integrable.
	We define the surface integral over $\Sigma$ as
	\begin{equation}
	\label{def:surf_int}
	\int_\Sigma\psi(\x )\,\dd S =\int_{\mathcal U}(\psi\circ\bs\Phi)(\bb u)\, |\bs{\tau}_1\times\bs{\tau}_2|\,\dd u^1\dd u^2 
	=\int_{\mathcal U}(\psi\circ\bs\Phi )(\bb u)\,\sqrt{g}\,\dd u^1\dd u^2,
	\end{equation}
	cp.\ Thomas et al.\ \cite[Section 16.6]{Thomas2004} or Bourne and Kendall \cite[Section 5.6]{Bourne1992}.
	
	Another important observation is that a vector may be decomposed into its tangential and normal components.
	This is due to the fact that the surface tangents and the surface normal form a basis of $\R^3$.
	For a vector field $\mathbf{a}: \Omega\to \R^3$ we thus may write
	\begin{align}
	    \mathbf{a} = a_\tau^\alpha\bs{\tau}_\alpha + a_\nu\bs{\nu}
			\label{eq:vector_decomp}
	\end{align}
	with the {\em tangential components} $a_\tau^\alpha$ and {\em normal component} $a_\nu =\bb a\cdot\bs\nu$. 
	
	When we move along the surface the tangents and the normal will change according to the coordinates. 
	First we want to discuss the tangents.
	From geometric intuition it is clear that in general the tangent will change in tangential and normal direction. 
	Thus we define according to (\ref{eq:vector_decomp}) a representation of the derivatives of $\bs{\tau}_\alpha$ as
	\begin{align}
	    \frac{\partial\bs{\tau}_\alpha}{\partial u_\beta} = \Gamma_{\alpha\beta}^{\gamma}\bs{\tau}_\gamma + b_{\alpha\beta}\bs{\nu}.
			\label{def:tangent_deriv}
	\end{align}
	The normal components $b_{\alpha\beta}$ of the tangent vectors are called the {\em curvature tensor} and are given as
	\begin{align}
	    b_{\alpha\beta} = \frac{\partial \bs{\tau}_\alpha}{\partial u_\beta}\cdot\bs{\nu} 
			= \frac{\partial^2\bs\Phi}{\partial u_\alpha\partial u_\beta}\bs{\nu}.
			\label{def:curv_tens}
	\end{align}
	Roughly speaking they give the change of the tangent vectors in normal direction with respect to the coordinates.
	Again the geometric view is quite clear since we expect a zero normal component of the derivative of a tangent for a flat surface,
	which has curvature zero.
	
	Along with the curvature tensor one can define the {\em mean curvature}\footnote{The factor $1/2$ coming from the
	arithmetic mean of the principal curvatures is the usual convention. It is also possible to define the mean curvature 
	without this factor, cf.\ Bothe \cite{Bothe2020,Bothe2022},
	Cermelli et al.\ \cite[Section 2]{Cermelli2005}, Grinfeld \cite[Section 12.4]{Grinfeld2013} 
	or Dziuk and Elliott \cite{Dziuk2007}.} $\kappa_M$ given by
	\begin{align}
	    \kappa_M = \frac{1}{2}g^{\alpha\beta}b_{\alpha\beta},
			\label{def:mean_curv}
	\end{align}
	see Aris \cite[(9.37.1),(9.38.7)]{Aris1989}.
	
	The tangential components in (\ref{def:tangent_deriv}) are called the {\em Christoffel symbols}. 
	They can be calculated by multiplying \eqref{def:tangent_deriv} with $g^{\delta\varepsilon}$ from the left
	and $\bs{\tau}_\delta$ from the right as well as using \eqref{def:surface_metric} to give
	\begin{align}
	\label{eq:chris}
	    g^{\delta\varepsilon}\frac{\partial\bs{\tau}_\alpha}{\partial u_\beta}\cdot\bs{\tau}_\delta 
			= \Gamma_{\alpha\beta}^{\gamma}g_{\gamma\delta}g^{\delta\varepsilon}
	    = \Gamma_{\alpha\beta}^\varepsilon.
	\end{align}
	We expand the left hand side and use the fact that $\bs{\tau}_\alpha$ was obtained as a derivative to obtain
  \begin{align*}			    
      g^{\delta\varepsilon}\frac{\partial\bs{\tau}_\alpha}{\partial u_\beta}\cdot\bs{\tau}_\delta
	    &=\frac{1}{2}g^{\delta\varepsilon}\left(\frac{\partial\bs{\tau}_\alpha}{\partial u_\beta}\cdot\bs{\tau}_\delta 
			+ \frac{\partial\bs{\tau}_\alpha}{\partial u_\beta}\cdot\bs{\tau}_\delta\right)
	    = \frac{1}{2}g^{\delta\varepsilon}\left(\frac{\partial\bs{\tau}_\alpha}{\partial u_\beta}\cdot\bs{\tau}_\delta
	    + \frac{\partial\bs{\tau}_\beta}{\partial u_\alpha}\cdot\bs{\tau}_\delta\right)\notag\\
	    &= \frac{1}{2}g^{\delta\varepsilon}\left(\frac{\partial(\bs{\tau}_\alpha\cdot\bs{\tau}_\delta)}{\partial u_\beta} 
			- \frac{\partial\bs{\tau}_\delta}{\partial u_\beta}\cdot\bs{\tau}_\alpha
	    + \frac{\partial(\bs{\tau}_\beta\cdot\bs{\tau}_\delta)}{\partial u_\alpha} 
			- \frac{\partial\bs{\tau}_\delta}{\partial u_\alpha}\cdot\bs{\tau}_\beta\right)\notag\\
	    &= \frac{1}{2}g^{\delta\varepsilon}\left(\frac{\partial g_{\alpha\delta}}{\partial u_\beta} 
			+ \frac{\partial g_{\beta\delta}}{\partial u_\alpha} -
	    \frac{\partial(\bs{\tau}_\alpha\cdot\bs{\tau}_\beta)}{\partial u_\delta}\right).
	\end{align*}
	With \eqref{eq:chris} we finally have, cp.\ Aris \cite[(9.91.2)]{Aris1989},
	\begin{align}
	    \Gamma_{\alpha\beta}^\gamma 
			&= \frac{1}{2}g^{\gamma\delta}\left(\frac{\partial g_{\alpha\delta}}{\partial u_\beta} 
			+ \frac{\partial g_{\beta\delta}}{\partial u_\alpha}
	    - \frac{\partial g_{\alpha\beta}}{\partial u_\delta}\right).
			\label{def:christoffel}
	\end{align}	

	The curvature tensor and the Christoffel symbols 
	are quantities which are intrinsic to the surface, see Aris \cite[Sections 9.13, 9.21]{Aris1989} 
	or Grinfeld \cite[Sections 5.9, 5.12, 10.9]{Grinfeld2013}.
	From \eqref{def:christoffel} or the fact that $\frac{\partial\bs\tau_\alpha}{\partial u_\beta} 
	=\frac{\partial\bs\tau_\beta}{\partial u_\alpha} 
	=\frac{\partial^2\bs\Phi}{\partial u_\alpha\partial u_\beta} $ in \eqref{def:tangent_deriv}
	it is clear that the Christoffel symbol is symmetric in the lower indices.
	
	Now we discuss how the change of the unit normal can be calculated.
	Intuitively there should be no change in their normal component. Indeed this can be seen differentiating $\bs{\nu}\cdot\bs{\nu} = 1$.
	We therefore may write
	\begin{align}
	    \frac{\partial\bs{\nu}}{\partial u_\alpha} = T^\beta_\alpha\bs{\tau}_\beta.
			\label{def:normal_derivative1}
	\end{align}

	In the following we want to determine the tangential components $T^\beta_\alpha$.
	Therefore we differentiate $\bs{\nu}\cdot\bs{\tau}_\gamma = 0$ with respect to $u^\alpha$
	\begin{align*}
	    0 = \frac{\partial\bs{\nu}}{\partial u_\alpha}\cdot\bs{\tau}_\gamma 
			+ \frac{\partial\bs{\tau}_\gamma}{\partial u_\alpha}\cdot\bs{\nu}
	    \quad\text{or equivalently}\quad
	    \frac{\partial\bs{\nu}}{\partial u_\alpha}\cdot\bs{\tau}_\gamma = -\frac{\partial\bs{\tau}_\gamma}{\partial u_\alpha}\cdot\bs{\nu}.
	\end{align*}
	If now (\ref{def:normal_derivative1}) is multiplied by $\bs{\tau}_\gamma$ and compared to the previous formula we obtain
	\begin{align*}
	    \frac{\partial\bs{\nu}}{\partial u_\alpha}\cdot\bs{\tau}_\gamma 
			= T^\beta_\alpha\bs{\tau}_\beta\cdot\bs{\tau}_\gamma = T^\beta_\alpha g_{\beta\gamma}
	    = -\frac{\partial\bs{\tau}_\gamma}{\partial u_\alpha}\cdot\bs{\nu}
	    \quad\text{or equivalently}\quad
	    T^\beta_\alpha = -g^{\beta\gamma}\frac{\partial\bs{\tau}_\gamma}{\partial u_\alpha}\cdot\bs{\nu}.
	\end{align*}
	Inserting this into (\ref{def:normal_derivative1}) and using (\ref{def:tangent_deriv}) finally gives
	\begin{align}
	    \frac{\partial\bs{\nu}}{\partial u_\alpha} = -g^{\beta\gamma}b_{\gamma\alpha}\bs{\tau}_\beta.
			\label{def:normal_derivative2}
	\end{align}
	\subsubsection*{Covariant Derivatives}

	So far we only used partial derivatives with respect to the surface coordinates. 
	This however might lead to problems since the result may depend on the coordinates chosen for some quantities.
	Therefore one needs a new derivative which is independent of the chosen parametrization. 
	This is called the {\em covariant surface derivative} $\nabla_\alpha$.
	We will omit the details of the derivation and interpretation. We again refer to 
	Aris \cite[Section 9.25]{Aris1989} and Grinfeld \cite[Chapter 11]{Grinfeld2013}.
	In short, one can think of the covariant derivative as respecting the change in the coordinates as well as the 
	change with respect to the tangents.
	
	In some generality the covariant derivative for a quantity with mixed indices having coordinates in the ambient space $\R^3$ 
	and on the surface we have, see Grinfeld \cite[(11.27)]{Grinfeld2013},
	\begin{align}
	    \nabla_\gamma T^{i\alpha}_{j\beta} = \frac{\partial T^{i\alpha}_{j\beta}}{\partial u^\gamma} 
			+ \bs{\tau}_{\gamma;k}\Gamma^i_{mk}T^{m\alpha}_{j\beta}
	    - \bs{\tau}_{\gamma;k}\Gamma^m_{jk}T^{i\alpha}_{m\beta} + \Gamma^\alpha_{\delta\gamma}T^{i\delta}_{j\beta} 
			- \Gamma^\delta_{\gamma\beta}T^{i\alpha}_{j\delta}.
			\label{def:covariant_deriv}
	\end{align}
	See also Aris \cite[Sections 9.25, 9.34]{Aris1989}.
	The basic rules are that for every upper index one has a $+$ Christoffel term and a $-$ 
	for every lower index. The term for every Latin index involves the tangent vector.
	If an index does not appear the whole Christoffel term with this index vanishes. For example the covariant derivative of
	a contravariant surface vector $\bb a$ or covariant spatial vector $\bb b$ are
	\begin{align*}
	    \nabla_\gamma a^\alpha = \frac{\partial a^\alpha}{\partial u^\gamma} + \Gamma^\alpha_{\delta\gamma}a^\delta
			\qquad\text{and}\qquad 	\nabla_\alpha b_j= \frac{\partial b_j}{\partial u_\alpha}-\tau_{\alpha;k}\Gamma^m_{jk}b_m.
	\end{align*}
	The covariant derivative of a scalar function $\psi$ is $\nabla_\alpha\psi =\frac{\partial \psi}{\partial u_\alpha}$.
	
	The covariant derivative has the following properties, more may be found in the given literature:
	\begin{enumerate}[(i)]
	    \item The covariant derivative coincides with the partial derivative when applied to 
			invariants, e.g.\ the normal $\bs{\nu}$ and its components.
			\item The covariant derivative of a tensor is again a tensor. The order of the tensor is increased by one.
	    \item The covariant derivative satisfies sum and product rules, but different derivatives do not commute.
	    \item The covariant derivative commutes with contraction.
	    \item The covariant derivatives $\nabla_\gamma g_{\alpha\beta}$ and $\nabla_\gamma g^{\alpha\beta}$ vanish. 
			This is also called \emph{metrinilic property}
			\footnote{Sometimes also called metrilinic property.}, see Grinfeld \cite[Subsection 8.6.7]{Grinfeld2013}.
	\end{enumerate}

	For the tangent vectors $\bs\tau_\alpha$ we obtain from \eqref{def:covariant_deriv} using \eqref{def:tangent_deriv}
	\begin{align}
	    \nabla_\beta\bs{\tau}_\alpha = \frac{\partial \bs{\tau}_\alpha}{\partial u_\beta} - \Gamma^\gamma_{\alpha\beta}\bs{\tau}_\gamma
	    = b_{\alpha\beta}\bs{\nu},
			\label{cova_deriv_tangent}
	\end{align}
	cp.\ Grinfeld \cite[(11.43)]{Grinfeld2013}.
	In this result the main property of the covariant derivative becomes apparent, 
	i.e.\ the result is independent from the chosen parametrization of the surface.
	Further we have that $\nabla_\beta\bs{\tau}_\alpha = \nabla_\alpha\bs{\tau}_\beta$.
	
	Differentiating $\bs\nu\cdot\bs\nu=1$ gives $\nabla_\alpha\bs\nu\cdot\bs\nu=0$. So $\nabla_\alpha\bs\nu$ is a tangent vector,
	i.e.\ $\nabla_\alpha\bs\nu =\nu^\beta_\alpha\bs\tau_\beta$. Further we have from $\bs\tau_\gamma\cdot\bs\nu=0$ that
	\[
	0=\nabla_\alpha \left(\bs\tau_\gamma\cdot\bs\nu\right)=\nabla_\alpha \bs\tau_\gamma\cdot\bs\nu+
	\bs\tau_\gamma\cdot\nabla_\alpha \bs\nu
	\]
	or using \eqref{cova_deriv_tangent}
	\[
	g_{\beta\gamma}\nu_\alpha^\beta=\nu_\alpha^\beta\bs\tau_\beta\cdot\bs\tau_\gamma =
	\nabla_\alpha\bs\nu\cdot\bs\tau_\gamma=-\nabla_\alpha\bs\tau_\gamma\cdot\bs\nu =-b_{\alpha\gamma}\bs\nu\cdot\bs\nu =-b_{\alpha\gamma}.
	\]
	This gives
	\begin{equation}
	\label{eq:deriv_nu}
	\nu_\alpha^\beta =-g^{\beta\gamma}b_{\gamma\alpha}\qquad\text{or}\qquad \nabla_\alpha\bs\nu
	= -g^{\beta\gamma}b_{\gamma\alpha}\bs\tau_\beta,
	\end{equation}
	cp.\ Aris \cite[(9.36.2),(9.36.3)]{Aris1989}.

	\subsection{Surface Differential Operators}
	\label{subsec:surf_diff}

	We now want to study surface differential operators.
	For this we assume that all mathematical pathologies are excluded, in particular that the 
	surface is sufficiently smooth. Let $\nabla_\x 
	= (\frac\partial{\partial\xin_1},\frac\partial{\partial\xin_2},\frac\partial{\partial\xin_3})=(\partial^1,\partial^2,\partial^3)$ 
	be the spatial gradient of a function.
	For a smooth function $\psi : \Omega\to \R$, with $\Omega\subseteq\R^3$ open,	
	we define the {\em surface gradient} as the projection\footnote{One can define the order 2, rank 1 
	tensor $\bb N=\left(\nu^i\nu_j\right)_{1\le i,j\le 3}$. It satisfies $\bb N^2 =\bb N$ since $\bs\nu\cdot\bs\nu =1$.
	Then for any covariant tensor $\bb a$ of order 1 we have
	$\bb N\bb a =(\bb a\cdot\bs\nu)\bs\nu=a_\nu\bs\nu$, i.e.\ $\bb N$ is the normal projection. With the
	Kronecker delta tensor $\bb I$ for $\R^3$ the complementary projection in the
	tangential direction is given by $\bb T =\bb I-\bb N$. This gives $\bb a^\shortparallel=\bb T\bb a$ and $\nabla_\Sigma =\bb T\nabla$.
	Mavrovouniotis and Brenner \cite[Section 2]{Mavrovouniotis1993} denoted $\bb T$ as the surface idemfactor.}
	in tangential direction
	\begin{align}
	    \begin{split}\label{def:surf_grad}
	        \nabla_\Sigma\psi &= \nabla_\x\psi - (\nabla_\x\psi\cdot\bs{\nu})\bs{\nu}
	    \end{split}
	\end{align}
	with components $\partial^i_\Sigma\psi = \partial^i\psi - (\nabla_\x\psi\cdot\bs{\nu})\nu^i$.
	For a surface quantity we have $\nabla_\Sigma\psi\cdot\bs{\nu} = 0$. 
	Setting $\psi =a_i$ in the $i$th component, we can define the {\em surface divergence} for 
	covariant vector fields $\bb a :\Omega \to \R^3$ as
	\begin{align}
	    \begin{split}\label{def:surf_div}
	        \nabla_\Sigma\cdot\bb a = \partial^i_\Sigma a_i
	        = \left[\partial^ia_i - (\nabla_\x a_i\cdot\bs{\nu})\nu^i\right]
	        = \nabla_\x\cdot\bb a - (\nabla_\x a_i\cdot\bs{\nu})\nu^i.
	    \end{split}
	\end{align}
	From both definitions we obtain, using the notations $\psi_\nu=\nabla_\x\psi\cdot\bs\nu$ and $a_\nu=\bb a\cdot\bs\nu$,
	the product rule
	\begin{align}
	\label{eq:product_rule}
	\nabla_\Sigma\cdot(\psi\bb a) =&\nabla_\x\cdot(\psi\bb a)-(\nabla_\x(\psi a_i)\cdot\bs\nu)\nu^i\notag\\
	=&\bb a\cdot\nabla_\x\psi+\psi\nabla_\x\cdot\bb a -(\nabla_\x\psi\cdot\bs\nu) a_i\nu^i-\psi (\nabla_\x a_i)\cdot\bs\nu)\nu^i\notag\\
	=&\bb a\cdot\nabla_\Sigma\psi +\psi_\nu a_\nu +\psi\nabla_\Sigma\cdot\bb a -\psi_\nu a_\nu 
	=\bb a\cdot\nabla_\Sigma\psi +\psi\nabla_\Sigma\cdot\bb a.
	\end{align}
	
	Now we want to deduce the equivalent coordinate formulation of the differential operators.
	Therefore, we calculate the tangential derivative using \eqref{def:surf_grad} and \eqref{def:surface_metric}
	\begin{align}
	\label{eq:grad}
		\nabla_\Sigma\psi \cdot\bs\tau_\gamma 
			&= \left(\nabla_\x \psi  - (\nabla_\x \psi \cdot\bs{\nu})\bs{\nu}\right)\cdot\bs\tau_\gamma
	    = \nabla_\x \psi \cdot\bs\tau_\gamma
			= \nabla_\x \psi \cdot\frac{\partial \bs\Phi}{\partial u^\gamma}\nonumber\\
			&= \frac{\partial \psi }{\partial u^\gamma}
			= \frac{\partial \psi }{\partial u_\alpha}\delta_\gamma^\alpha 
			= g^{\alpha\beta}\frac{\partial \psi }{\partial u_\alpha}g_{\beta\gamma} 
			= g^{\alpha\beta}\frac{\partial \psi }{\partial u_\alpha}
			\bs\tau_\beta\cdot\bs\tau_\gamma.
  \end{align}
	Note that we found the natural relations
	\begin{equation}
	\label{eq:grad_rel}
	\nabla_\Sigma\psi \cdot\bs\tau_\alpha =\frac{\partial \psi }{\partial u_\alpha}=\nabla_\alpha\psi \qquad\text{and}\qquad
	\nabla_\Sigma\psi \cdot\bs\tau_\alpha =\nabla_\x \psi \cdot\bs\tau_\alpha.
	\end{equation}
	We also have $\nabla_\Sigma\psi \cdot\bs\nu= 0$ and all three vectors $\bs\tau_1,\bs\tau_2,\bs\nu$
	form a basis of $\R^3$. Therefore, \eqref{eq:grad} implies that
	\begin{align}
			\label{eq:surface_gradient}
	    \nabla_\Sigma\psi	= g^{\alpha\beta}\frac{\partial \psi}{\partial u_\alpha}\bs\tau_\beta
			= g^{\alpha\beta}\frac{\partial \psi}{\partial u_\alpha}\frac{\partial\bs\Phi}{\partial u_\beta}
	\end{align}
	with components 
	%
	$\partial_\Sigma^i\psi 
	= g^{\alpha\beta}\frac{\partial \psi}{\partial u_\alpha}\frac{\partial\Phi^i}{\partial u_\beta}
	= g^{\alpha\beta}\frac{\partial \psi}{\partial u_\alpha}\tau_{\beta;i}
	$. 
	Thus we can again take $\psi=a_i$ and write the surface divergence of a covariant vector field $\mathbf{a}$ as
	\begin{equation}
	\label{eq:surf_div}
	    \nabla_{\Sigma}\cdot\bb a=\partial_\Sigma^ia_i(\mathbf{x}) 
			= g^{\alpha\beta}\frac{\partial a_i}{\partial u_\alpha}\frac{\partial\Phi^i}{\partial u_\beta}
			= g^{\alpha\beta}\frac{\partial a_i}{\partial u_\alpha}\tau_{\beta;i}.
	\end{equation}
	This can be compared to the surface divergence defined in Aris \cite[(9.41.6)]{Aris1989}. The latter can be written as
	\begin{align}
	    \nabla_\Sigma\cdot \mathbf{a} = g^{\alpha\beta}\bs{\tau}_\beta\nabla_\alpha\mathbf{a}
			= g^{\alpha\beta}\frac{\partial\Phi^i}{\partial u_\beta}\nabla_\alpha a_i.
			\label{def:surf_divergence_spcevec}
	\end{align}
	Since the covariant derivative of a scalar function is $\nabla_\alpha a_i =\frac{\partial a_i}{\partial u_\alpha}$ both
	\eqref{eq:surf_div} and \eqref{def:surf_divergence_spcevec} are equivalent. 
		
	So far we considered derivatives of functions $\psi$ and vector fields $\bb v$ defined in the volume surrounding the surface. 
	The formulas \eqref{def:surf_grad} and \eqref{def:surf_div}
	are only defined for such functions and vector fields respectively. 
	But we can also consider functions and vector fields that are only defined on the surface
	$\Sigma$, i.e.\ $\psi:\Sigma\to\R$. We then say that such a function is smooth iff $\psi(\bs \Phi(\cdot)):{\mathcal U}\to\R$ is a
	smooth function. For such functions the surface gradient \eqref{eq:surface_gradient} is well defined. So thereby the surface gradient 
	\eqref{def:surf_grad} may be extended to such functions. We can proceed analogously for vector fields 
	and define their divergence via \eqref{eq:surf_div}.
	
	For further calculations it is important to obtain the surface divergence of  the normal $\bs{\nu}$. 
	We start with \eqref{def:surf_divergence_spcevec}, use \eqref{eq:deriv_nu}, \eqref{def:surface_metric} 
	and \eqref{def:mean_curv} to obtain
	\begin{align}
	    \nabla_\Sigma\cdot\bs{\nu} &= g^{\alpha\beta}\bs{\tau}_\alpha\cdot\nabla_\beta\bs{\nu} 
			= -g^{\alpha\beta}\bs{\tau}_\alpha\cdot g^{\varepsilon\gamma}b_{\gamma\beta}\bs\tau_\varepsilon
	    = -g^{\alpha\beta}(\bs{\tau}_\alpha\cdot\bs{\tau}_\varepsilon)g^{\varepsilon\gamma}b_{\gamma\beta}\notag\\
	    &= -g^{\alpha\beta}g_{\alpha\varepsilon}g^{\varepsilon\gamma}b_{\gamma\beta}
	    = -g^{\alpha\beta}\delta_\alpha^\gamma b_{\gamma\beta}
	    = -g^{\alpha\beta}b_{\alpha\beta}= -2\kappa_M.
			\label{eq:weingarten_normal}
	\end{align}
		
	Let us decompose the vector field $\bb a:\R^3\to\R^3$ into a {\em tangential part} 
	$\bb a^\shortparallel$ and a {\em normal part} $\bb a^\perp$ as
	$\bb a = \bb a^\shortparallel + \bb a^\perp$, cp.\  \eqref{eq:vector_decomp}, with
	\begin{equation}
	\label{eq:decomp}
			\bb a^\perp = (\bb a\cdot\bs{\nu})\,\bs{\nu} =a_\nu\,\bs\nu ,\qquad\text{and}\qquad
	    \bb a^\shortparallel = \bb a - \bb a^\perp =a_\tau^\alpha\bs\tau_\alpha.
	\end{equation}
	This implies $\bb a^\shortparallel\cdot\bs{\nu} = 0$. 
	We can reformulate the surface divergence as follows. We calculate the divergence of the normal part using \eqref{def:surf_div}
	\begin{equation}
	\label{eq:perp_div}
	\nabla_\Sigma\cdot\bb a^\perp =\nabla_\Sigma\cdot\left[(\bb a\cdot\bs{\nu})\,\bs{\nu}\right]
	= \nabla_\Sigma(\bb a\cdot\bs{\nu})\cdot\bs{\nu} +(\bb a\cdot\bs{\nu})\nabla_\Sigma\cdot\bs{\nu}.
	\end{equation}		
	With \eqref{eq:surface_gradient} we determine the first term to be
	\begin{equation}
	\label{eq:vanish}
	\nabla_\Sigma(\bb a\cdot\bs{\nu})\cdot\bs{\nu}
	=g^{\alpha\beta}\frac{\partial (a_j\nu_j)}{\partial u_\alpha}\bs\tau_\beta\cdot\bs\nu =0.
	\end{equation}
	So now \eqref{eq:perp_div} gives $\nabla_\Sigma\cdot\bb a^\perp =(\bb a\cdot\bs{\nu})\nabla_\Sigma\cdot\bs{\nu}$. 
	With this and \eqref{eq:weingarten_normal} we obtain
	\begin{align}
	    \nabla_\Sigma\cdot\bb a &= \nabla_\Sigma\cdot\left(\bb a^\shortparallel + \bb a^\perp\right)
	    = \nabla_\Sigma\cdot\bb a^\shortparallel + (\bb a\cdot\bs{\nu})\nabla_\Sigma\cdot\bs{\nu} 
			= \nabla_\Sigma\cdot\bb a^\shortparallel -2 \kappa_M a_\nu ,
			\label{def:surf_divergence_spacevec}
	\end{align}
	cp.\ Dziuk and Elliot \cite[(2.6)]{Dziuk2007}.
		
	In the following we want to derive a further useful identity for the surface divergence of the tangential part of a vector field.
	Let us first consider the product $\bs{\nu}\cdot (\nabla_\x\times(\bs{\nu}\times\bb a))$. We apply the 
	$\epsilon-\delta$-identity \eqref{eq:ep_delta}, $\bs\nu\cdot\bs\nu=1$, $\nabla_\x(\bs\nu\cdot\bs\nu)=0$ 
	and \eqref{def:surf_div} to obtain
	\begin{align*}
	    \bs{\nu}\cdot (\nabla_\x\times(\bs{\nu}\times\bb a)) &= \nu^i\epsilon_{ijk}\partial^j(\epsilon_{krs}\nu_r a_s)
			=\nu^i\epsilon_{kij}\epsilon_{krs}\partial^j(\nu_r a_s)\\
			&=\nu^i (\delta_{ir}\delta_{js}-\delta_{is}\delta_{jr})\partial^j(\nu_r a_s)
			= \nu^i\partial^j(\nu^i a_j)-\nu_i\partial^j(\nu_j a_i) \\
			&=\nu^ia_j\partial^j\nu_i +\nu^i\nu_i \partial^ja_j-\nu^ia_i\partial^j\nu_j -\nu^i\nu_j \partial^ja_i\\
			&= \bb a\cdot\frac 12\nabla_\x (\bs\nu\cdot\bs\nu)  +  \nabla_\x\cdot\bs a
			-(\bb a\cdot\bs\nu)\nabla_\x\cdot\bs\nu-(\nabla_\x a_i\cdot\bs\nu)\nu^i\\ 
			&=\nabla_\Sigma\cdot\bb a-(\bb a\cdot\bs\nu)\nabla_\x\cdot\bs\nu.
	\end{align*}
	Now we compute the surface divergence of the tangential part of the vector field $\bb a$. We use
	the definition \eqref{eq:decomp} and the identity \eqref{eq:vanish} to get
	\begin{align*}
	    \nabla_\Sigma\cdot\bb a^\shortparallel &= \nabla_\Sigma\cdot\left(\bb a - (\bb a\cdot\bs{\nu})\bs{\nu}\right)
	    = \nabla_\Sigma\cdot\bb a 
			-\nabla_\Sigma(\bb a\cdot\bs{\nu})\cdot\bs{\nu}- (\bb a\cdot\bs{\nu})\nabla_\Sigma\cdot\bs{\nu}\\
	    %
	    %
	    %
	    &= \nabla_\Sigma\cdot\bb a - (\bb a\cdot\bs{\nu})\nabla_\x\cdot\bs{\nu}.
	\end{align*}
	Combining the two calculations, we obtain the formula
	\begin{align}
	    \nabla_\Sigma\cdot\bb a^\shortparallel = \left(\nabla_\x\times(\bs{\nu}\times\bb a)\right)\cdot\bs{\nu}.
			\label{eq:surf_div_tangent}
	\end{align}
	\subsection{Moving Surfaces}
	\label{subsec:move}

	So far we only considered a fixed surface $\Sigma \subset \Omega$ with $\Omega\subseteq\R^3$ open. 
	Now we want to extend the results to moving surfaces.
	A moving surface $\Sigma$ can be considered as a family of surfaces $\Sigma(\t )$ with a parameter $\t  \in [0, T]=I$
	for some $T > 0$. The stationary {\em initial reference surface} $\Sigma_0 = \Sigma(0)$ may be described as before via a smooth 
	bijective parametrization $\bs\Phi_{\Sigma_0} : \mathcal{U} \to \Sigma_0$ with $\mathcal U\subseteq\R^2$ an open subset.
	We further assume that we have a smooth mapping 
	$\bs\chi_\Sigma:I\times\Sigma_0\to \R^3$ with $\bs\chi_\Sigma(\t,\y)\in \Sigma(\t)$. 
	We define the {\em flow map} $\bs\chi_\Sigma^\t(\y) =\bs\chi_\Sigma(\t,\y)$ of the surface and assume that $\bs\chi_\Sigma^\t$ 
	is a diffeomorphism, i.e.\ a bijective smooth map, from
	the initial surface $\Sigma_0$ to $\Sigma(\t)$. 
	So it parametrizes the moving surface $\Sigma(\t)$ over the stationary initial surface $\Sigma_0$.
	Note that $\bs\chi_\Sigma^0=\bs\chi_\Sigma(0,\cdot )$ is the identity map
	on $\Sigma_0$ and that $\bs\chi_{\Sigma}$ is not directly related to the map $\bs\chi$ 
	considered in Subsection \ref{subsec:control}. The latter represents infinitesimal particles that 
	may move through the surface $\Sigma(\t)$. But $\bs\chi_{\Sigma}$ describes the movement of points on the surface.
			 
	Now we may also describe the surface $\Sigma(\t )$ in terms of the surface parameters 
	$\bb u=(u^1,u^2) ^T\in \mathcal{U}$ using the mapping
	\begin{align*}
	    \bs\Phi = \bs\chi_{\Sigma}(\cdot ,\bs\Phi_{\Sigma_0}(\cdot)) : I\times\mathcal{U} \to \Sigma(\t ).
	\end{align*}
	The results from the previous subsections remain valid and we assume that the dependence on the parameter $t$ is smooth.
	We introduce the {\em velocity} 
	of a {\em surface point} $\mathbf{x}\in\Sigma(\t)$ that is given by
	\begin{align}
	    \bb w(\t,\x)=\frac{\partial\bs\chi_\Sigma}{\partial\t}(\t,\y)
			\qquad\text{for}\quad\x =\bs\chi_\Sigma(\t,\y)\in \Sigma(\t)
			\label{def:surf_velocity}
	\end{align}
	or
	\begin{equation}
	\label{eq:surf_velocity}
			\mathbf{w}(\t ,\bs\Phi(\t,\bb u)) = \frac{\partial\chi_{\Sigma}}{\partial \t }(\t ,\bs\Phi_{\Sigma_0}(\bb u))
			= \frac{\partial\bs\Phi}{\partial \t }(\t ,\bb u).
	\end{equation}
	It describes the movement of the surface points and can be independent of the surrounding particle flow velocity $\bb v$
	that will be introduced in Subsection \ref{subsec:control}. 	
	We have for the trajectory in $\R^3$ of such a point on the surface
	\[
	  \frac{\dd\mathbf{x}}{\dd \t }(\t) = \mathbf{w}(\t ,\mathbf{x}(\t)) 
		= \mathbf{w}(\t ,\bs\Phi(\t ,\bb u)) .
	\]

	The parametrization and the flow map are mathematical tools used in order to have a better understanding of the surface
	movement and to prove theorems. We only need to have a parametrization of the initial surface and the flow map.
	In practical applications
	the velocity or the normal velocity of the surface will be determined from physical relations. The moving surface geometry
	or the flow map would have to be determined by integrating such vector fields. Mathematical theory 
	for this integration was given by Bothe \cite{Bothe2020}.
	
	As in \eqref{eq:decomp} we can decompose the velocity as $\mathbf{w} = \mathbf{w}^\shortparallel + \mathbf{w}^\perp$.
	The vector $\mathbf{w}^\shortparallel$ is sometimes called {\em tangential coordinate velocity}.
	The term $w_\nu =\bb w\cdot\bs\nu$ can be called {\em velocity of the surface} since the normal direction is implied.
	Then $w_\nu\bs{\nu}$ is the {\em vector normal velocity}, see Grinfeld \cite[(15.32)]{Grinfeld2013},
	and $w_\nu$ the normal speed.
	It is important to note that $w_\nu$ is an invariant of the surface. A nice description and geometric motivation is given in 
	Grinfeld \cite[Sections 15.4, 15.5]{Grinfeld2013}.

	\subsubsection*{Implicit Description and an Explicit Parametrization of Surfaces}

	Let us assume that a surface $\Sigma (\t)\subset\Omega$ is given as an implicit function $f(\t,\x)=0$ with 
	$f:I\times  \R^3\to\R$ and $|\nabla_\x f|>0$, which implies $|\nabla_{\t,\x}f|>0$. 
	The time and space unit normal vector field is given as $\bb n = \frac{\nabla_{\t,\x}f}{|\nabla_{\t,\x}f|}$
	since the gradient is orthogonal to level sets. We have from \eqref{def:surf_velocity}
	\[
	0 =\frac\dd{\dd\t} f(\t,\bs\chi_\Sigma (\t,\x)) 
	= \left(f_\t +\nabla_\x f\cdot \frac{\partial\bs\chi_\Sigma}{\partial\t}\right)(\t,\bs\chi_\Sigma (\t,\x))
	= \left(f_\t +\nabla_\x f\cdot \bb w\right)(\t,\bs\chi_\Sigma (\t,\x)).
	\]
	and for $i=1,2,3$
	\[
	0 =\frac{\partial}{\partial u_\alpha} f(\t,\bs\Phi (\t,\bb u)) 
	= \left(\nabla_\x f\cdot \frac{\partial\bs\Phi}{\partial u_\alpha}\right)(\t,\bs\Phi (\t,\bb u))
	= \left(\nabla_\x f\cdot \bs\tau_\alpha\right)(\t,\bs\Phi (\t,\bb u)).
	\]
	The latter equation proves that $\bs\nu = \frac{\nabla_\x f}{|\nabla_\x f|}$ is 
	the spatial unit normal vector field that is orthogonal to the tangent vectors $\bs\tau_\alpha$ in $\R^3$ 
	for any $\t\in I$. The former gives with the notation $f_t=\frac\partial{\partial\t}f$
	\begin{equation}
	\label{eq:normal_vel}
	w_\nu =\bb w\cdot \bs\nu =\bb w\cdot \frac{\nabla_\x f}{|\nabla_\x f|}=\frac{-f_t}{{|\nabla_\x f|}}.
	\end{equation}

	Suppose that we can describe the points $\x\in\Sigma(\t)$ on the surface by a function $\sigma:I\times\Omega^2\to\R$
	with a suitable set $\Omega^2\subseteq\R^2$ as $\x=(\xin_1,\xin_2,\sigma(\t,\xin_1,\xin_2))$.
	Then it is given implicitly by $0=f(\t,\x)=\xin_3-\sigma(\t,\xin_1,\xin_2)$
	for $(\xin_1,\xin_2)\in\Omega^2$. We set $\x_\sigma =(\xin_1,\xin_2,\sigma(\t,\xin_1,\xin_2))\in\R^3$. 
	For the normal velocity $w_\nu$ of the surface $\Sigma$ we have from \eqref{eq:normal_vel} 
	\[
	w_\nu = \frac{-f_\t}{|\nabla_\x f|}= \frac{\sigma_\t}{\sqrt{(\sigma_{\xin_1})^2+(\sigma_{\xin_2})^2+1}}.
	\]
	This implies that
	\begin{equation}
	\label{eq:sigma_t}
	\sigma_\t(\t,\xin_1,\xin_2) =w_\nu (\t,\x_\sigma)\sqrt{(\sigma_{\xin_1})^2+(\sigma_{\xin_2})^2+1}.
	\end{equation}
	\subsubsection*{Intrinsic Time Derivatives}

	A major issue when working with moving surfaces is how to define an intrinsic time derivative that is invariant under
	a change of parametrization. We cannot keep $\x$ fixed, since $\x$ has moved when we consider
	a time $\t+\Delta\t$. For this purpose
	Thomas \cite[Section 4]{Thomas1957}, \cite[Section II,3]{Thomas1961} introduced a time derivative 
	for functions $\psi:I\times\Sigma(\t)\to\R$ 
	defined at $\x\in\Sigma(\t)$ as
	follows. Consider the unit normal vector $\bs\nu(\x)$ and the straight line 
	\[
	\mathcal L=\{\z\in\R^3 |\z=\lambda \bs\nu_\Sigma(\x)+\x \;\text{for}\;\lambda\in\R\}.
	\]
	For a small time $\Delta\t$ we take the moved surface $\Sigma(\t+\Delta\t)$ and let $\z^\x$ be the point where $\mathcal L$ intersects
	$\Sigma(\t+\Delta\t)$. Then we define the {\em Thomas time derivative}, assuming that the limit exists,
	\begin{equation}
	\label{def:thomas_deriv}
	\frac \delta{\delta\t} \psi (\t,\x)= \lim_{\Delta\t\to 0}\frac{\psi(\t+\Delta\t,\z^\x)-\psi(\t,\x)}{\Delta\t}.
	\end{equation}
	Truesdell and Toupin \cite[Section 179]{Truesdell1960}, Bowen and Wang \cite{Bowen1971}, Petryk and Mroz \cite{Petryk1986} 
	as well as Bothe \cite{Bothe2022} defined this time derivative by taking the curve that is
	obtained by integrating the normal velocity field $w_\nu\bs\nu_\Sigma$ instead of the line $\mathcal L$.
	Both curves have the velocity $w_\nu(\t,\x)\bs\nu_\Sigma(\t,\x)$ at $\x$. So both approaches give the same result.
	The definition is independent of a parametrization of the surface. Truesdell and Toupin \cite[Section 179]{Truesdell1960}
	as well as Bowen and Wang \cite{Bowen1971} discussed the extension of this derivative to tensors.
	
	We can interpret this derivative by assuming that $\widehat{\psi}:I\times\widehat{\Omega}\to\R$
	is an extended smooth function 
	with $\Sigma(t)\subset\widehat{\Omega}\subseteq\Omega\subseteq\R^3$ an open neighborhood of $\Sigma(t)$. 
	An simple extension of $\psi$ to a neighborhood 
	of the surface may be obtained by a function that is constant in the directions normal to the surface
	on both sides of the surface, see e.g.\ Cermelli et al.\ \cite[Subsection 3.4]{Cermelli2005}. 
	This would lead to a vanishing normal derivative. Thus, we then have 
	$\mathbf{w}\cdot\nabla_\x\widehat{\psi}_\Sigma = w^\shortparallel\cdot\nabla_\Sigma\widehat{\psi}_\Sigma$ 
	and $\widehat{\psi}_\nu =0$. We do not assume such a specific extension, so we may have $\widehat{\psi}_\nu\ne 0$.
	Let $\bb u_{\Delta\t}$ be chosen so that $\z^\x =\bs\Phi(\t+\Delta\t,\bb u_{\Delta\t})$.	
	Inserting the parametrization $\x=\bs\Phi(\t,\bb u)$ into \eqref{def:thomas_deriv} implies that
	\begin{align}
	\label{eq:thomas_deriv}
	\frac \delta{\delta\t} \widehat{\psi} (\t,\bs\Phi(\t,\bb u))
		=&\lim_{\Delta\t\to 0}\frac{\widehat{\psi}(\t+\Delta\t,\bs\Phi(\t+\Delta\t,\bb u_{\Delta\t}))
	-\widehat{\psi}(\t,\bs\Phi(\t,\bb u))}{\Delta\t}\\
	=&\lim_{\Delta\t\to 0}\frac{\widehat{\psi}(\t+\Delta\t,\bs\Phi(\t+\Delta\t,\bb u_{\Delta\t}))
	-\widehat{\psi}(\t+\Delta\t,\bs\Phi(\t,\bb u))}{\Delta\t}\notag\\
	&+\lim_{\Delta\t\to 0}\frac{\widehat{\psi}(\t+\Delta\t,\bs\Phi(\t,\bb u))
	-\widehat{\psi}(\t,\bs\Phi(\t,\bb u))}{\Delta\t}\notag\\
	=&\lim_{\Delta\t\to 0}\bigg(\sum_{i=1}^3\frac{\widehat{\psi}(\t+\Delta\t,\bs\Phi(\t+\Delta\t,\bb u_{\Delta\t}))
	-\widehat{\psi}(\t+\Delta\t,\bs\Phi(\t,\bb u))}{\Phi^i(\t+\Delta\t,\bb u_{\Delta\t})-\Phi^i(\t,\bb u)}\notag\\
	&\qquad\qquad\qquad\cdot\frac{\Phi^i(\t+\Delta\t,\bb u_{\Delta\t})-\Phi^i(\t,\bb u)}{\Delta\t}\bigg)
	+ \frac\partial{\partial\t}\widehat{\psi}(\t,\bs\Phi(\t,\bb u))\notag\\
	=&\left[\nabla_\x\widehat{\psi}\cdot w_\nu\bs\nu_\Sigma + 	\frac\partial{\partial\t}\widehat{\psi}\right](\t,\bs\Phi(\t,\bb u))
	=\left[\frac\partial{\partial\t}\widehat{\psi}+\widehat{\psi}_\nu w_\nu	\right](\t,\bs\Phi(\t,\bb u)).\notag
	\end{align}
	We see that for functions that are only defined on a surface the Thomas derivative
	is a generalization of the partial derivative with respect to time for such extended functions,
	as the term $\widehat{\psi}_\nu$ will vanish if the extension by constants along normal directions is used.
			
	For $\psi:I\times\Sigma\to\R$ we can proceed analogously to the Thomas time derivative
	when we replace the normal velocity $w_\nu(\t,\x)\bs\nu_\Sigma(\t,\x)$ by the
	velocity $\bb w(\t,\x)$ of the surface. We can use the flow map and 
	set $\widehat{\z}^\x =\bs\chi^{\Delta\t}_\Sigma (\x)$.
	Then we define 	the {\em Lagrangian time derivative} for surfaces, assuming that the limit exists,
	\begin{equation}
	\label{def:lagrange_deriv}
	\frac D{D\t} \psi (\t,\x)= \lim_{\Delta\t\to 0}\frac{\psi(\t+\Delta\t,\widehat{\z}^\x)-\psi(\t,\x)}{\Delta\t}.
	\end{equation}
	With $\Delta\z^\x=\widehat{\z}^\x-\z^\x$ we then obtain
	\begin{align}
	\label{eq:lagrange_deriv}
	\frac D{D\t} \psi (\t,\x)=& \lim_{\Delta\t\to 0}\frac{\psi(\t+\Delta\t,\widehat{\z}^\x)-\psi(\t+\Delta\t,\z^\x)}{\Delta\t}
	+ \lim_{\Delta\t\to 0}\frac{\psi(\t+\Delta\t,\z^\x)-\psi(\t,\x)}{\Delta\t}\notag\\
	=& \lim_{\Delta\t\to 0}\sum_{i=1}^3\frac{\psi(\t+\Delta\t,\widehat{\z}^\x)-\psi(\t+\Delta\t,\z^\x)}{\Delta\zin_i^\x}
	\frac{\Delta\zin_i^\x}{\Delta\t} + \frac \delta{\delta\t} \psi (\t,\x)\notag\\
	=& \frac \delta{\delta\t} \psi (\t,\x)+\bb w^\shortparallel (\t,\x)\cdot\nabla_\Sigma \psi(\t,\x).
	\end{align}
	So the Lagrangian derivative $\frac D{D\t}$ is the Thomas derivative $\frac \delta{\delta\t}$ plus an advective term
	\footnote{In Grinfeld \cite[(15.34)]{Grinfeld2013} the Lagrangian time derivative is written as a partial time derivative. The Thomas
	time derivative is introduced as the invariant time derivative $\dot{\nabla}$. Therefore, in the definition
	of $\dot{\nabla}$ via \cite[(15.37)]{Grinfeld2013} the advective term has a minus sign.}.
	
	Intrinsic time derivatives for general tensors can be found in Truesdell and Toupin \cite[Section 179]{Truesdell1960}.
	They are called displacement derivatives there. For a more recent discussion of these derivatives
	see Bowen and Wang \cite{Bowen1971}.
	
	Let us set $\widetilde{\psi}(\t,\bb u) =\psi(\t,\bs\Phi (\t,\bb u))$ 
	for a function $\psi:I\times\Sigma\to\R$. We can define, assuming that the limit exists, the time derivative
	\begin{equation}
	\label{def:mueller_deriv}
	\mathring{\psi}(\t,\bb u)=\frac\partial{\partial\t}\widetilde{\psi}(\t,\bb u)
	=\lim_{\Delta\t\to 0}\frac{\widetilde{\psi}(\t+\Delta\t,\bb u)-\widetilde{\psi}(\t,\bb u)}{\Delta\t}
	\end{equation}
	This derivative was used by Dreyer \cite[(48)]{Dreyer2003} and M\"uller \cite[(3.9)]{Mueller1985} in transport theorems for surfaces.
	Using the definition of $\widetilde{\psi}$ and $\widehat{\z}^\x =\bs\chi^{\Delta\t}_\Sigma (\bs\Phi (\t,\bb u))
	=\bs\Phi (\t+\Delta\t,\bb u)$ it satisfies
	\begin{align}
	\label{eq:mueller_deriv}
	\mathring{\psi}(\t,\bb u)=&\frac\partial{\partial\t}\widetilde{\psi}(\t,\bb u)
	=\lim_{\Delta\t\to 0}\frac{\psi(\t+\Delta\t,\bs\Phi (\t+\Delta\t,\bb u))-\psi(\t,\bs\Phi (\t,\bb u))}{\Delta\t}\notag\\
	=&\lim_{\Delta\t\to 0}\frac{\psi(\t+\Delta\t,\widehat{\z}^\x)-\psi(\t,\x)}{\Delta\t}=\frac D{D\t}\psi (\t,\x).
	\end{align}
	So the derivative \eqref{def:mueller_deriv} is equivalent to the Lagrangian time derivative \eqref{def:lagrange_deriv}.
	This justifies the notation $\mathring{\psi}(\t,\x)=\frac D{D\t}\psi (\t,\x)$ for the Lagrangian time derivative.

	\subsubsection*{The Time Derivative of the Metric}

	Since the metric is time dependent now we also want to calculate the 
	time derivative of the determinant $g$ of the metric. We apply \eqref{eq:mat_det_deriv} to the symmetric matrix
	$\mathbf{G}(\t ) = g_{\alpha\beta}(t)$ with $\mathbf{G}(\t )^{-1} = g^{\alpha\beta}(t)$ and $\det(\mathbf{G}) = g$. This gives,
	using the symmetry of $\bb G(\t )$,
	\begin{equation}
	\label{eq:g_deriv}
	\frac{\dd g}{\dd \t }  = g^{\alpha\beta}\,\frac{\dd g_{\alpha\beta}}{\dd \t }\, g.
	\end{equation}
	By using \eqref{def:surface_metric} and then \eqref{def:tangent_normal} we have
	\[
	\frac{\dd g_{\alpha\beta}}{\dd \t }=\frac{\partial^2\bs\Phi}{\partial u_\alpha\partial\t}
	\frac{\partial\bs\Phi}{\partial u_\beta}+\frac{\partial\bs\Phi}{\partial u_\alpha}
	\frac{\partial^2\bs\Phi}{\partial u_\beta\partial\t} 
	= \frac{\partial w_k}{\partial u_\alpha} \frac{\partial\Phi^k}{\partial u_\beta}
	+\frac{\partial\Phi^k}{\partial u_\alpha} \frac{\partial w_k}{\partial u_\beta}.
		\]
	With \eqref{eq:g_deriv}, \eqref{eq:surf_div} and the symmetry of the inverse metric tensor $g^{\alpha\beta}$ 
	we have the following formula
	\begin{align}
	\label{eq:metricdet_deriv}
	    \frac{\dd}{\dd \t }\sqrt{g} =& \frac{1}{2\sqrt{g}}\,\frac{\dd}{\dd \t }g
			= \frac{1}{2\sqrt{g}}\,g^{\alpha\beta}\,\frac{\dd g_{\alpha\beta}}{\dd \t }\, g 
			= \frac 12\left(g^{\alpha\beta}\,\frac{\partial w_k}{\partial u_\alpha}	\frac{\partial\Phi^k}{\partial u_\beta} 
			+ g^{\beta\alpha}\frac{\partial w_k}{\partial u_\beta}\frac{\partial\Phi^k}{\partial u_\alpha}\right)\,\sqrt g
			\nonumber\\
	    =& (\nabla_\Sigma\cdot\mathbf{w})\,\sqrt{g}.
	\end{align}
	This is the time derivative of the surface element needed in the next section when the integrals \eqref{def:surf_int} 
	are taken over moving surfaces and differentiated with respect to time. It is quite natural that the evolution
	of this element is determined by the surface divergence of the velocity field of the surface.
	
	Bothe \cite{Bothe2020,Bothe2022} considered moving interfaces in a stationary domain. In \cite{Bothe2020} he gave various
	results for the trajectories 
	$\x(\t)=\bs\chi^\x(\t)=\bs\chi(\t,\y)$ for $\t\in I$
	of points on the surface given as solutions to systems of ordinary differential equations.

	\subsection{Moving Control Volumes with Singular Moving Surfaces}
	\label{subsec:control}
	\subsubsection*{Motions and Control Volumes}

For the description of particle flows in $\R^n$ one may consider a smooth map
$\bs{\chi }:I \times \R ^n \rightarrow \R ^n$ with $\bs \chi  (0,\x ) = \x \in\R^n$. 
By smooth we mean that the map is continuously differentiable as often 
as needed. We assume it to be at least a diffeomorphism. For any fixed $\y\in \R ^n$ the map defines a trajectory
$\bs\chi ^{\x}:I\to \R^n , $ via $\x(\t)=\bs \chi^{\x} (\t)=\bs\chi (\t,\y)$. 
These are the trajectories of infinitesimal particles that are initially
at $\y$ and then at $\x(\t)$ for $\t\in [0,T]$. No particles should appear or disappear. 
So we assume that the map $\bs{\chi }^{\t }:\R ^n \rightarrow \R ^n$, defined as $\bs{\chi }^{\t }(\y )
=\bs{\chi }(\t ,\y ) $ for any fixed $\t  \in I $ must be a smooth bijection.
We have $\bs{\chi }^{0}(\y )=\y $, i.e.\ this map is the identity map in $\R^n$.
We define
\[
\bb{v}(\t ,\x )= \frac{\partial}{\partial \t}  \bs\chi (\t ,\x ) 
\]
to be the {\em velocity field} of the particles.
This description of a particle flow is known as the formulation in {\em Lagrangian} coordinates $\bs\chi $.
In this section we assume that it is continuous. The mappings $\bs\chi$ or the family of mappings $\bs\chi^\t$ 
are called {\em motions} or {\em flow maps}.
The invertible Jacobian matrices $\bb D_\x\bs\chi^\t=\left(\frac{\partial\bs\chi^\t_j}{\partial\xin_k}\right)_{1\le j,k\le 3}$ 
are the {\em deformation gradients}.
		
Let ${\mathcal V} (0)\subset\R^n$ be any open subset.	We call it the {\em reference volume}	
and ${\mathcal V} (\t ) = \bs{\chi }^{\t }[{\mathcal V} (0)]\subset\R^n$ a {\em material control volume}.
This means that we assume that there is no mass flux through the boundary 
$\mathcal S (\t )=\partial\mathcal{V}(\t )$ of $\mathcal{V}(\t )$.
Thus the boundary $\mathcal S (\t)$ moves with the velocity of the particles sitting on it.
With particle one does not actually refer to the molecules or atoms of the material which occupy a positive volume.
In a continuum at any time $\t$ we interpret any spatial point as the location of an {\em infinitesimal particle}.
The dynamics of the control volume and its boundary are given by the space time Lagrangian mapping $\bs\chi$. 
For simplicity we will not always emphasize the time dependence of such a Lagrangian volume $\mathcal{V}$ or surface $\mathcal{S}$.
A stationary Eulerian volume will later be denoted by $\Omega \subseteq \R^n$.	
	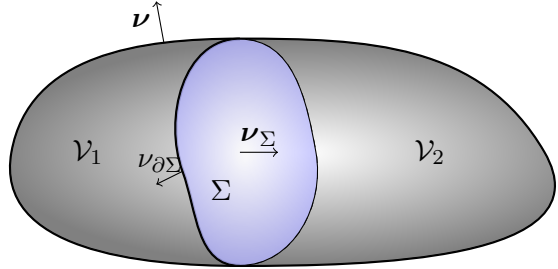
\begin{figure}[!ht]
	    \centering
	    \begin{tikzpicture}
	        %
	        \draw[very thick] (0,0) to[out=80, in=180] (3,1.5) to[out=0,in=120] (7,0) 
					to[out=-60,in=0] (3,-1.5) to[out=180,in=-100] (0,0) -- cycle;
	        \shade[inner color = white, outer color=gray, draw=black] (0,0) to[out=80, in=180] (3,1.5) 
					to[out=0,in=120] (7,0) to[out=-60,in=0] (3,-1.5) to[out=180,in=-100] (0,0) -- cycle;
	        %
	        \draw (1,0) node {$\mathcal{V}_1$};
	        \draw (5.5,0) node {$\mathcal{V}_2$};
	        \draw[->] (2,1.45) -- (1.9,2);
	        \draw (1.7,1.75) node {$\bs\nu$};
	        %
	        \draw[dashed] (3,-1.5) to[out=0,in=-80] (4,0) to[out=100,in=0] (3,1.5);
	        \draw[very thick] (3,-1.5) to[out=180,in=-70] (2.25,-0.25) to[out=110,in=180] (3,1.5);
	        \shade[inner color = white, outer color=blue!30, draw=black, fill opacity=0.8] (3,-1.5) 
					to[out=0,in=-80] (4,0) to[out=100,in=0] (3,1.5) to[in=110,out=180] (2.25,-0.25) to[in=180,out=-70] (3,-1.5) -- cycle;
	        %
	        \draw (2.75,-0.5) node {$\Sigma$};
	        \draw [->] (3,0) -- (3.5,0);
	        \draw (3.25,0.2) node {$\bs{\nu}_\Sigma$};
	        \draw [->] (2.25,-0.25) -- (1.9,-0.425);
	        \draw (1.95,-0.15) node {$\nu_{\partial\Sigma}$};
	    \end{tikzpicture}
	    \caption{An arbitrary material volume $\mathcal{V}$ separated by an immaterial internal surface $\Sigma$ into
	    the sub-volumes $\mathcal{V}_1$ and $\mathcal V_2$. }
	    \label{fig:gen_volume}
	\end{figure}
	\subsubsection*{Moving Surfaces in Control Volumes}

	Let us consider such a control volume $\mathcal V\subset\R^3$
	that contains a singular moving surface $\Sigma$ separating it 
	into two open volumes denoted by $\mathcal{V}_1$ and $\mathcal{V}_2$.
	By singular we mean that physical states may be discontinuous across the surface.
	We have $\mathcal V=\mathcal V_1\cup\Sigma\cup\mathcal V_2$.
	Such a situation is sketched in Figure \ref{fig:gen_volume}.
	Note that whereas the boundary surface $\mathcal S$ is assumed to be a material surface, the internal surface $\Sigma$ 
	may allow particle transfer between the two regions. We assume that we have a smooth mapping 
	$\bs\chi_\Sigma:I\times\Sigma_0\to \R^3$ from the initial reference surface $\Sigma_0=\Sigma(0)$ to the
	moving surface $\Sigma=\Sigma(\t)$, as introduced in Subsection \ref{subsec:move}.
	The velocity field $\bb w$ on $\Sigma$ is independent of the velocity field of particles defined above.
	
	We take $\bs\nu$ to be outwards pointing unit normal vector field on $\mathcal S=\partial\mathcal{V}$.
	The unit normal field on the surface $\Sigma$ pointing from $\mathcal{V}_1$ to $\mathcal{V}_2$ is denoted by $\bs{\nu}_\Sigma$.
	Let $\bs\tau$ be a tangent vector field to the boundary curve
	$\partial\Sigma$. We assume that the unit normal vector field $\bs\nu_\Sigma$ to the surface and $\bs\tau$ satisfy 
	the usual right hand rule. Looking at the surface from the side to which the normal points, 
	the curve is oriented in the counter-clockwise direction. 
	We set $\bs\nu_{\partial\Sigma} = \bs{\tau}\times\bs{\nu}_\Sigma$ to be the outward pointing boundary normal 
	that is tangential to the surface.
	
	With out loss of generality we assume that all components of the normal field $\bs\nu_\Sigma$ are positive
	and that $w_\nu>\bb v\cdot\bs\nu$. So the volume $\mathcal V_1$ becomes larger
	at the expense of $\mathcal V_2$ since the surface $\Sigma$ is moving faster in the normal direction 
	than the local material velocity of particles and into  $\mathcal V_2(\t)$. 
	We denote by $\widehat{\mathcal V}_i(\t)=\bs\chi^t[\mathcal V_i(0)]$ for $i=1,2$ the material volumes of the particles 
	that are initially in $\mathcal V_i(0)$ and introduce the difference volume
	$\mathcal V_D(\t)=\mathcal V_1(\t)\backslash\widehat{\mathcal V}_1(\t)$. So the volume $\mathcal V_1(\t)$ is divided into two parts.
	The initial set $\mathcal V_2(0)$ is also divided into two parts
	$\mathcal V_{21}^\t\subset\mathcal V_2(0)$ and $\mathcal V_{22}^\t\subset\mathcal V_2(0)$. 
	Note that they contain infinitesimal particles at time $0$. But their composition changes with time.
	The set $\mathcal V_{21}^\t$ 
	contains the particles that move from $V_2(0)$ to $\mathcal V_1\cup\Sigma$ in the time interval $[0,\t]$
	i.e.\ $\mathcal V_D(\t)=\bs\chi^\t[\mathcal V_{21}^\t]$.
	In $\mathcal V_{22}^\t$ are the particles that make up $\mathcal V_2(\t) =\bs\chi^\t[\mathcal V_{22}^\t]$.
	
	Now we want to parametrize the moving boundary $\Sigma_0(\t)$ between $\mathcal V_{21}$ and $\mathcal V_{22}$.
	We have $\Sigma(\t)=\bs\chi^t[\Sigma_0(\t)]$ and assume that its unit normal vector field  $\bs\nu_{\Sigma_0}$ points towards 
	$\mathcal V_{22}^\t$.
	We denote by $\bs\chi^{-\t}$ the inverse of the bijective map $\bs\chi^\t$.
	For $\t\in I$ let $\bs\chi_{\Sigma_0}(\t,\cdot)=\bs\chi^{-\t}(\t,\bs\chi_\Sigma(\t,\cdot)):\Sigma_0\to \R^3$.
	Then for $\y\in\Sigma_0$ we have $\bs\chi_{\Sigma_0}(\t,\y)\in \Sigma_0(\t)$ 
	and $\bs\chi_\Sigma(\t,\y) =\bs\chi(\t,\bs\chi_{\Sigma_0}(\t,\y))$. The
	mapping $\bs\chi_{\Sigma_0}$ parametrizes the moving boundary $\Sigma_0(\t)$ between $\mathcal V_{21}^\t$ and $\mathcal V_{22}^\t$
	over the stationary surface $\Sigma_0$.
	
	With the invertible Jacobian matrix or deformation gradient $\bb D_\x\,\bs\chi^\t  
	=\left( \frac{\partial\bs{\chi }^{\t }_{j}}{\partial x _{k}}\right)_{1\leq j,\, k \leq n }$ we have
	\[
	\frac{\partial}{\partial\t}\bs\chi_\Sigma(\t,\y) =\frac{\dd}{\dd\t}\bs\chi(\t,\bs\chi_{\Sigma_0}(\t,\y))
	=\frac{\partial\bs\chi}{\partial\t}(\t,\bs\chi_{\Sigma_0}(\t,\y))
	+\bb D_\x\,\bs\chi^\t(\t,\bs\chi_{\Sigma_0}(\t,\y))\frac{\partial\bs\chi_{\Sigma_0}}{\partial\t}(\t,\y).
	\]
	Let $\y=\bs\chi_{\Sigma_0}(\t,\y_0)$, $\x=\bs\chi^\t(\y)$ and $\bb w_0(\t,\y)=\frac{\partial\bs\chi_{\Sigma_0}}{\partial\t}(\t,\y_0)$.
	Also we have $\bb w(\x)=\frac{\partial}{\partial\t}\bs\chi_\Sigma(\t,\y_0)$ and
	$\bb v(\x)=\frac{\partial\bs\chi}{\partial\t}(\t,\bs\chi_{\Sigma_0}(\t,\y_0))$.
	Then we obtain
	\[
	\bb w(\t,\x) =\bb v(\t,\x)	+\bb D_\x\,\bs\chi^\t(\t,\y)\bb w_0(\t,\y).
	\]
	Solving for $\bb w_0(\t,\y)$ we have 
	\begin{equation}
	\label{eq:u_def}
	\bb w_0(\t,\y)=\left(\bb D_\x\,\bs\chi^\t\right)^{-1}(\t,\x)\left[\bb w(\t,\x)
	-\bb v(\t,\x)\right].
	\end{equation}
	The velocity of $\chi_{\Sigma_0}$ is a transformation of the relative velocity $\bb w-\bb v$ of the surface $\Sigma(\t)$.

	\subsubsection*{Integration Formulas}

	We want to assume that the boundary $\mathcal S_2(0)$ of the initial control volume $\mathcal V_2(0)$ can be parametrized over 
	the closure $\overline{\Omega}\subseteq\R^2$ of a fixed open shadow domain $\Omega$ in the following way.
	There are suitable functions $\alpha,\beta:\Omega\to \R$ with $\alpha(\yin_1,\yin_2)=\beta(\yin_1,\yin_2)$ 
	for $(\yin_1,\yin_2)\in\partial\Omega$, $\alpha(\yin_1,\yin_2)<\beta(\yin_1,\yin_2)$ for $(\yin_1,\yin_2)\in\Omega$ and
	\begin{equation}
	\label{eq:v2_para}
	\mathcal V_2(0) =\{\,\y\in\R^3\, |\, (\yin_1,\yin_2)\in\Omega\;\text{and}\;\alpha(\yin_1,\yin_2)<\yin_3<\beta(\yin_1,\yin_2)\,\}.
	\end{equation}
	The functions $\alpha$, $\beta$ parametrize the boundary $\mathcal S_2(0)$. 
	Under these circumstances we have for any integrable function $h:\mathcal V_2(0)\to\R$ the formula.
	\begin{equation}
	\label{eq:int_para}
	\int_{\mathcal V_2(0)} h(\y)\,\dd V_0 =\int_\Omega\int_{\alpha(\yin_1,\yin_2)}^{\beta(\yin_1,\yin_2)}h(\y)\dd \yin_3\dd A.
	\end{equation}
	For many arguments with integrals over control volumes
	one can make such an assumption. One can possibly use another selection of the two variables or divide the control
	volume into pieces that satisfy such an assumption.
	
	We already assume that $\nu^3_\Sigma>0$ and we additionally assume that $\nu^3_{\Sigma_0}> 0$. 
	This implies that the moving surface $\Sigma_0(\t)$ can be described by a function $\sigma:I\times\Omega\to\R$.
	Since we define it over $\Omega$ there may be parts of $\Sigma_0(\t)$ that lie outside of $\mathcal V_2(0)$. 
	But this does not matter.
	The variables $\yin_1$, $\yin_2$ correspond to the parameters $u^1$, $u^2$ in Subsection \ref{subsec:surf}.
	We have the parametrization $\bs\Phi_{\Sigma_0}(\yin_1,\yin_2) =(\yin_1,\yin_2,\sigma(\t,\yin_1,\yin_2))$ with $\t$ as an
	extra parameter. The tangent vectors are $\bs\tau_1^0=(1,0,\sigma_{\yin_1})$ and $\bs\tau_2^0=(0,1,\sigma_{\yin_2})$
	giving
	\[
	\bs\tau_1^0\times\bs\tau_2^0 = (-\sigma_{\yin_1},-\sigma_{\yin_2},1)\qquad\text{and}\qquad
	\bs\nu_{\Sigma_0} =\frac{\bs\tau_1^0\times\bs\tau_2^0}{|\bs\tau_1^0\times\bs\tau_2^0|}
	=\frac{(-\sigma_{\yin_1},-\sigma_{\yin_2},1)}{\sqrt{(\sigma_{\yin_1})^2+(\sigma_{\yin_2})^2+1}}.
	\]
	The surface integrals \eqref{def:surf_int} become, using the plane surface element $\dd A=\dd\yin_1\dd\yin_2$, 
	\begin{equation}
	\label{eq:sigma_0_int}
	\int_{\Sigma_0} \psi(\x)\dd S =\int_\Omega \psi(\bs\Phi_{\Sigma_0}(\y))|\bs\tau_1^0\times\bs\tau_2^0|\,\dd A
	=\int_\Omega \psi(\bs\Phi_{\Sigma_0}(\y))\sqrt{(\sigma_{\yin_1})^2+(\sigma_{\yin_2})^2+1}\,\dd A.
	\end{equation}
	
	We obtain a parametrization of the singular surface $\Sigma(\t)$ by the mapping $\bs\Phi_\Sigma=\bs\chi^\t\circ\bs\Phi_{\Sigma_0}$.
	It gives the tangent vectors as $\bs\tau_\alpha =\bb D_\y\bs\chi^\t\bs\tau_\alpha^0$ for $\alpha=1,2$.
	We use the determinant $\mathcal J^\t$ of the Jacobian matrix of $\bs\chi^\t$ as given by \eqref{def:J}. 
	From equation \eqref{eq:prod_trans} we have 
	\[
	\bs\tau_1\times \bs\tau_2 =\mathcal J^\t(\bb D_\x\bs\chi^{-\t})^T\left(\bs\tau_1^0\times \bs\tau_2^0\right)
	\]
	and using \eqref{def:tangent_normal} this gives
	\[
	\bs\nu_\Sigma|\bs\tau_1\times \bs\tau_2 |=\bs\tau_1\times \bs\tau_2 
	=\left((\bb D_\x\bs\chi^{-\t})^T\bs\nu_{\Sigma_0}\right)\,\mathcal J^\t\,|\bs\tau_1^0\times \bs\tau_2^0|.
	\]
	For any vector field $\bb a$ with $\bb a\cdot\bs\nu_\Sigma$ integrable on $\Sigma(\t)$ we obtain the integral transformation
	\begin{equation}
	\label{eq:int_surf_trans}
	\int_{\Sigma} \bb a\cdot\bs\nu_\Sigma\, \dd S=
	\int_\Omega\bb a\cdot\bs\nu_\Sigma|\bs\tau_1\times \bs\tau_2 |\,\dd A
	=\int_\Omega\left((\bb D_\x\bs\chi^{-\t})\bb a\right)\cdot\bs\nu_{\Sigma_0}\mathcal J^\t|\bs\tau_1^0\times \bs\tau_2^0|\,\dd A.
	\end{equation}
	\section{Generic Balance Laws}
	\label{sec:balance_law}

	In this short section we want to consider generic balance laws that may be used to describe the 
	dynamics of a fluid in the presence of interfaces.
	A general balance equation for an additive, i.e.\ extensive, physical quantity $\Psi$ 
	states that the total change in time is equal to the flux across the boundary and the sources
	\begin{align}
	\label{eq:gen_balance_law}
	    \tikzmarkin[fill=green,opacity=0.2]{eq0:box1a}(0,-0.3)(0,0.6)
	    \tikzmarkin[fill=cyan,opacity=0.2]{eq0:box1b}(0,-0.7)(0,0.3)
	    \frac{\dd\Psi}{\dd \t }  
	    \tikzmarkend{eq0:box1b}
	    \tikzmarkend{eq0:box1a}
			=
			\tikzmarkin[fill=green,opacity=0.2]{eq0:box2a}(0,-0.3)(0,0.6)
	    \tikzmarkin[fill=cyan,opacity=0.2]{eq0:box2b}(0,-0.7)(0,0.3)
			\underbrace{\Phi}_{\text{Flux}} 
			\tikzmarkend{eq0:box2b}
	    \tikzmarkend{eq0:box2a}
			+
			\tikzmarkin[fill=green,opacity=0.2]{eq0:box3a}(0,-0.3)(0,0.6)
	    \tikzmarkin[fill=cyan,opacity=0.2]{eq0:box3b}(0,-0.7)(0,0.3)
			\underbrace{\Xi}_{\text{Source}}
			\tikzmarkend{eq0:box3b}
	    \tikzmarkend{eq0:box3a}.
	\end{align}
	We take a control volume $\mathcal V\subset\R^3$ that is separated
	into $\mathcal V_1$ and $\mathcal V_2$ by a surface $\Sigma$, as in Figure \ref{fig:gen_volume}.
	On $\mathcal V$ we have a volume density $\psi$ and we assume an additional
	surface density $\psi_\Sigma$ on $\Sigma$. Then $\Psi$ may be written as
	\begin{align}
	    \Psi(\t ) =
	    \tikzmarkin[fill=green,opacity=0.2]{eq1:box1a}(0,-0.3)(0,0.6)
	    \tikzmarkin[fill=cyan,opacity=0.2]{eq1:box1b}(0,-0.7)(0,0.3)
	    \int\limits_{\mathcal{V}_1\cup\mathcal{V}_2}\psi\,\dd V
	    \tikzmarkend{eq1:box1a}
	    \tikzmarkend{eq1:box1b}
	    + 
			\tikzmarkin[fill=cyan,opacity=0.2]{eq1:box2}(0,-0.7)(0,0.6)
	    \int\limits_\Sigma \psi_\Sigma\,\dd S
	    \tikzmarkend{eq1:box2}.
			\label{psi_decomp}
	\end{align}
	We use the color green to mark volume quantities and blue for surface quantities. The terms marked with mixed
	colors are volume quantities that will later contribute to the surface balance.
	
	We introduce a {\em volume flux density} $\bb j$ and a {\em surface 
	flux density} $\bb j_\Sigma$. But, we write them with boundary flux integrals as
	\begin{align}
	    \Phi(\t ) =
	    -\tikzmarkin[fill=green,opacity=0.2]{eq2:box1a}(0,-0.3)(0,0.6)
	    \tikzmarkin[fill=cyan,opacity=0.2]{eq2:box1b}(0,-0.7)(0,0.3)
	    \int\limits_{\mathcal{S}} \bb j\cdot\bs\nu\,\dd S
	    \tikzmarkend{eq2:box1a}
	    \tikzmarkend{eq2:box1b}
	    -\tikzmarkin[fill=cyan,opacity=0.2]{eq2:box2}(0,-0.7)(0,0.6)
	    \int\limits_{\partial\Sigma} \bb j_\Sigma\cdot\bs\nu_{\partial\Sigma}\,\dd l
	    \tikzmarkend{eq2:box2},
			\label{phi_decomp}
	\end{align}
	as in Dreyer \cite[Subsection 5.1]{Dreyer2003} and M\"uller \cite[Section 3.1]{Mueller1985}.
	The fluxes may include differential terms such as Fick's or Fourier's laws. The resulting differential equations
	can therefore be diffusive, i.e.\ parabolic.
	
	The first integral accounts for the non-convective flux across the closed surface of any 
	arbitrarily small material volume.
	Whereas, the second integral describes the flux which is tangential to the surface $\Sigma$ 
	and normal to its closed boundary $\partial\Sigma$. The unit normal vector field $\bs\nu_{\partial\Sigma}$ 
	on the curve 	$\partial\Sigma$ is chosen
	to be tangential to the surface $\Sigma$ and outward with respect to the surface $\Sigma$.
	Since the unit normal vector fields $\bs\nu$ and $\bs\nu_{\partial\Sigma}$ are pointing outwards, 
	the sign takes care that incoming fluxes contribute positively to $\Phi$. Using the Gauss-Green theorem
	both boundary integrals can later be converted to appropriate domain integrals of the divergence of the fluxes.
		
	Finally, the source $\Xi$ can be decomposed using a {\em volume density} $\xi$ and a {\em surface density} 
	$\xi_\Sigma$. Thus it is given as
	\begin{align}
	    \Xi(\t ) =
	    \tikzmarkin[fill=green,opacity=0.2]{eq3:box1}(0,-0.7)(0,0.6)
	    \int\limits_{\mathcal{V}_1\cup\mathcal{V}_2}\xi\,\dd V
	    \tikzmarkend{eq3:box1}
	    +\tikzmarkin[fill=cyan,opacity=0.2]{eq3:box2}(0,-0.7)(0,0.6)
	    \int\limits_\Sigma \xi_\Sigma\,\dd S
	    \tikzmarkend{eq3:box2}.
			\label{xi_decomp}
	\end{align}

	In the case of a zero right hand side, the equation (\ref{eq:gen_balance_law}) becomes a physical conservation law in a strict
	sense. This is the case for conservation of mass leading to the continuity equation.
	For other generic quantities $\Psi$ we have the additional terms on the right hand side making them balance laws.
	The boundary flux density integrals with $\bb j$ and $\bb j_\Sigma$ can be transformed into divergence integrals over the volume
	that together with the terms from the left hand side  of \eqref{eq:gen_balance_law} lead to partial differential equations 
	in divergence form, then it has become a custom that the resulting differential equations are called conservation laws. Such 
	differential equations are called balance laws if they additionally
	contain non-differential source terms or differential terms that cannot be brought into divergence form.
	
	Inserting equations (\ref{psi_decomp}), (\ref{phi_decomp}) and (\ref{xi_decomp}) into the 
	general balance law (\ref{eq:gen_balance_law}) gives
	\begin{align}
	    \frac{\dd}{\dd \t }\left(
	    \tikzmarkin[fill=green,opacity=0.2]{eq4:box1a}(0,-0.3)(0,0.6)
	    \tikzmarkin[fill=cyan,opacity=0.2]{eq4:box1b}(0,-0.7)(0,0.3)
	    \int\limits_{\mathcal{V}_1\cup\mathcal{V}_2}\psi\,\dd V
	    \tikzmarkend{eq4:box1a}
	    \tikzmarkend{eq4:box1b}\right.
	    &+ \left. \tikzmarkin[fill=cyan,opacity=0.2]{eq4:box2}(0,-0.7)(0,0.6)
	    \int\limits_\Sigma \psi_\Sigma\,\dd S
	    \tikzmarkend{eq4:box2}
	    \right)\notag\\
	    &= -\tikzmarkin[fill=green,opacity=0.2]{eq4:box3a}(0,-0.3)(0,0.6)
	    \tikzmarkin[fill=cyan,opacity=0.2]{eq4:box3b}(0,-0.7)(0,0.3)
	    \int\limits_{\mathcal{S}} \bb j\cdot\bs\nu\,\dd S
	    \tikzmarkend{eq4:box3a}
	    \tikzmarkend{eq4:box3b}
	    -\tikzmarkin[fill=cyan,opacity=0.2]{eq4:box4}(0,-0.7)(0,0.6)
	    \int\limits_{\partial\Sigma} \bb j_\Sigma\cdot\bs\nu_{\partial\Sigma}\,\dd l
	    \tikzmarkend{eq4:box4}
	    + \tikzmarkin[fill=green,opacity=0.2]{eq4:box5}(0,-0.7)(0,0.6)
	    \int\limits_{\mathcal{V}_1\cup\mathcal{V}_2}\xi\,\dd V
	    \tikzmarkend{eq4:box5}
	    +\tikzmarkin[fill=cyan,opacity=0.2]{eq4:box6}(0,-0.7)(0,0.6)
	    \int\limits_\Sigma \xi_\Sigma\,\dd S
	    \tikzmarkend{eq4:box6}.
			\label{eq:int_balance_law}
	\end{align}
	Examples for physical quantities and their corresponding densities may be found in Bedeaux \cite{Bedeaux1986},
	Dreyer \cite[Section 6]{Dreyer2003} as well as Truesdell and Toupin \cite{Truesdell1960}. Some are nicely summarized 
	in a table in M\"uller \cite[Section 3.2]{Mueller1985}.
	
	The basic assumption that surface quantities contribute to the total quantity $\Psi$ is a major point in our considerations.
	In fact, this assumption will be responsible for additional terms in the equations on singular surfaces
	with discontinuous weak solutions. This includes differential terms, see e.g.\ (\ref{eq:gen_loc_jump_cond}), 
	that are absent in the standard case of shock or contact discontinuities where the jump conditions are purely algebraic.

	\section{Transport Theorems, Volume and Surface Balance Equations}
	\label{sec:transport}

	In the following the time derivative on the left hand side of equation (\ref{eq:int_balance_law}) 
	is evaluated using transport theorems since the coordinates of volume and surface depend on time.
	The resulting relations then lead to differential equations in the volume and on the surface.

	\subsection{Volume Transport}
	\label{subsec:vol_trans}

	We start with the well known Reynolds transport theorem for volume integrals,
	see also Aris \cite[Section 4.22]{Aris1989}, Truesdell and Toupin \cite[Section 81]{Truesdell1960} 
	as well as Warnecke \cite[Satz 1.3]{Warnecke1999}. The proof is given since it allows some insight into the result
	and its later extensions. We state the theorem for any dimension $n=1,2,3,\ldots$ since it does not depend on dimension.
	For dimension $n=1$ \eqref{eq:reynolds_thm1} gives versions of Leibniz's rule for integrals, 
	see e.g.\ Thomas \cite[Chapter 5]{Thomas2004} for a version
	that is missing the integral term arising due to a parameter dependent integrand.
	\begin{thm}[The Reynolds Transport Theorem in $\R^n$]
	\label{thm:reynolds_transport}
	    Let $\psi : I\times\mathbb{R}^n \to \mathbb{R}$ be a continuously differentiable density function
	    and $\mathcal{V}(\t ) \in \mathbb{R}^n$ be a closed volume moving with 
			the local particle velocity $\mathbf{v}(\t ,\mathbf{x})$. We assume that the boundary $\mathcal S(\t)=\partial\mathcal V(\t)$
			consists of finitely many smooth pieces. Let $\bs\nu$ denote the outer unit normal field of the
			hypersurface $\mathcal S(\t)$. Then the following equations hold
	    \begin{equation}
			\label{eq:reynolds_thm}
	        %
	        \frac{\dd}{\dd \t }\Psi(\t)
					= 
					\frac{\dd}{\dd \t }\int\limits_{\mathcal{V}(\t )}\psi(\t,\mathbf{x})\,\dd V
	        =
					\int\limits_{\mathcal{V}(\t )}\left[\frac{\partial}{\partial \t }\psi
					+ \nabla_\x\cdot(\psi\mathbf{v})\right](\t ,\x) \,\dd V
      %
			\end{equation}
			and
			\begin{equation}
			\label{eq:reynolds_thm1}
					\frac{\dd}{\dd \t }\Psi(\t)
					=
					\int\limits_{\mathcal{V}(\t )}\frac{\partial}{\partial \t }\psi(\t ,\mathbf{x}) \,\dd V
					+ 
					\int\limits_{\mathcal{S}(\t )}\left[\psi\,\mathbf{v}\cdot\bs\nu\right](\t ,\x)\,\dd S
			\end{equation}		
	\end{thm}
	\begin{prf}
	    The key idea is that at any time the volume $\mathcal{V}(\t)$ may be described applying 
			a smooth transformation to the initial volume $\mathcal{V}(0)$.
	    Let the initial volume $\mathcal{V}(0)$ be described using the coordinates $\y = (\yin_1, \dots, \yin_n)^T$.
	    Since we consider a closed material volume no particles are removed from or added to the volume.
	    Thus we assume the existence of a diffeomorphism $\bs\chi:I\times\R^n$ such that at every time $\t\in I$ 
			we have $\x(\t) = \bs\chi(\t ,\y_0)=\bs\chi^\t(\y_0)$, see Appendix \ref{app:det}. 
			This is the Lagrangian description of the flow or motion. 
			
	    The Jacobian matrix or deformation gradient of $\bs\chi^\t$ for fixed $\t $ is given by $\textbf{\textup{D}}\bs\chi^\t(\y) 
			= (\partial\bs\chi^\t_j/\partial \yin_k)_{1\le j,k\le n}$. Further, let
	    ${\mathcal J}^\t = \det(\textbf{\textup{D}}\bs\chi)$ denote the corresponding determinant.
	    The determinant ${\mathcal J}^\t$ is needed to relate the volume element $\dd V$ at time $\t$ to the volume element 
			of the initial configuration $\dd V_0$, i.e.\ $\dd V = {\mathcal J}^\t\dd V_0$.
	    The derivative with respect to time of  ${\mathcal J}^\t$ is given by \eqref{eq:det_deriv} as
	    $\frac{\dd {\mathcal J}^\t}{\dd \t } = (\nabla_\x\cdot\mathbf{v})\,{\mathcal J}^\t$ Now we obtain
	    \begin{align}
			\label{eq:reynolds_proof}
	        \frac{\dd}{\dd \t }\Psi(\t ) &= \frac{\dd}{\dd \t }\int\limits_{\mathcal{V}(\t )}\psi(\t ,\y)\,\dd V
	        = \frac{\dd}{\dd \t }\int\limits_{\mathcal{V}(0)}\psi(\t ,\bs\chi^\t(\y)){\mathcal J}^\t(\t ,\y)\,\dd V_0\notag\\
	        &= \int\limits_{\mathcal{V}(0)}\frac{\dd}{\dd \t }\left[\psi(\t ,\bs\chi^\t(\y)){\mathcal J}^\t(\t,\y)\right]
					\,\dd V_0\notag\\
	        &= \int\limits_{\mathcal{V}(0)}\left[\frac{\dd}{\dd \t }\psi(\t ,\bs\chi^\t(\y)) 
					+ \psi(\t ,\bs\chi^\t(\y))[\nabla_\x\cdot\mathbf{v}(\t ,\bs\chi^\t(\y))]\right]{\mathcal J}^\t(\t ,\y)\,\dd V_0\notag\\
	        &= \int\limits_{\mathcal{V}(\t )}\left[\frac{\dd}{\dd \t }\psi(\t ,\mathbf{x}) 
					+ \psi(\t ,\mathbf{x})(\nabla_\x\cdot\mathbf{v}(\t ,\mathbf{x}))\right]\,\dd V.
	    \end{align}
	    Using the chain rule and $\mathbf{v}(\t ,\y) = \dd\y/\dd \t  = \partial\bs\chi^\t (\y)/\partial \t $ leads to
	    \begin{align*}
	        \frac{\dd}{\dd \t }\psi(\t ,\x ) &= \frac{\partial}{\partial \t }\psi(\t ,\mathbf{x}) 
					+ \nabla_\x\psi(\t ,\mathbf{x})\cdot\mathbf{v}(\t ,\mathbf{x}).
	    \end{align*}
			This inserted into the last integral gives \eqref{eq:reynolds_thm}. The second equation \eqref{eq:reynolds_thm1} follows
			by applying the Gauss-Green Theorem to the integral with the divergence for any arbitratry but fixed $\t$.
			In this derivation $\t$ is just a parameter.
	\end{prf}

	Note that the divergence term in \eqref{eq:reynolds_thm} is a combination of two contributions, of the advective 
	or Lagrangian derivative 
	resulting from the total time derivative of $\psi$ and of the divergence of the velocity field coming from the time derivative of
	the volume element. The first part is the advection of the 
	physical quantity $\psi$ along stream lines of the flow. The second is from the local expansion or compression of the control volume
	due to the flow. Both combine nicely to give the divergence term which is mathematically very appreciated when weak solutions 
	are considered in conjunction with jump discontinuities of physical quantities in flows, 
	such as shock waves, contact discontinuities, reaction fronts 
	or phase transitions.

	\subsection{Divergence Theorem and Transport Theorem for Embedded Moving Surfaces}
	\label{subsec:surf_transport}

	The case $n=2$ of the Reynolds Theorem \ref{thm:reynolds_transport} would only cover moving surfaces in the plane $\R^2$.
	In order to treat moving surfaces inside of moving control volumes in $\R^3$ we need an additional transport theorem.
	First, we give a divergence or Gauss-Green theorem for curved surfaces in $\R^3$.
	For this we consider any open subset $\Omega\subseteq\R^3$ and a bounded smooth surface $\Sigma \subset \R^3$ that is 
	bounded by a smooth curve $\partial\Sigma$.
	Let $\bs\tau$ be a tangent vector field to the boundary curve.
	We assume that the vector field $\bs\nu_\Sigma$, which is normal to the surface, and the tangent 
	vector field $\bs\tau$ satisfy the usual right hand rule. We 
	take $\bs\nu_{\partial\Sigma} = \bs{\tau}\times\bs{\nu}_\Sigma$ to be the outward pointing boundary normal 
	that is tangential to the surface. Now we want to show the divergence theorem for surfaces in $\Omega$, cp.\
	Lee \cite[Theorem 16.32]{Lee2013} for a more general version of \eqref{eq:surface_div}.
	\begin{thm}[Divergence Theorem or Gauss-Green Theorem for Surfaces in $\R^3$]
	\label{thm:surface_int}
	        Let $\Sigma \subset \Omega\subseteq\R^3$ be a bounded smooth surface.
	        Further, $\bb a:\Sigma\to\R^3$ is a continuously differentiable vector field that is either defined on the
					boundary $\partial\Sigma$ or has a bounded continuous extension to this boundary.
	        Like in \eqref{eq:decomp} it may be decomposed into tangential and normal components
					as follows $\bb a = \bb a^\shortparallel + a_\nu\bs\nu_\Sigma$. By $\dd l$ we denote the line element on 
					the curve $\partial \Sigma$. We assume that the curve is continuous and consists of finitely many
					smooth pieces.
	        Then the following divergence formula for surface integrals holds
	        \begin{align}
	            \int\limits_\Sigma \left[\nabla_\Sigma\cdot\bb a^\shortparallel\right](\x)\;\dd S
							= \int\limits_{\partial\Sigma} \left[\bb a\cdot\bs\nu_{\partial\Sigma}\right](\x)\,\dd l .
	            \label{eq:surface_div}
	        \end{align}
					From this we obtain the formula
	        \begin{align}
	            \int\limits_\Sigma \left[\nabla_\Sigma\cdot\bb a\right](\x)\;\dd S
							= \int\limits_{\partial\Sigma} \left[\bb a\cdot\bs\nu_{\partial\Sigma}\right](\x)\,\dd l 
							-\int\limits_\Sigma\left[ 2\kappa_Ma_\nu\right](\x)\;\dd S.
	            \label{eq:surface_div_2}
	        \end{align}
	\end{thm}
	\begin{prf}
		  The proof is basically an application of the Kelvin-Stokes theorem, cp.\
			Aris \cite[(3.32.2)]{Aris1989}, Bourne \cite[Sections 6.5,6.6]{Bourne1992}, Grinfeld\cite[(14.53)]{Grinfeld2013}
			or Thomas et al.\ \cite[Section 16.7]{Thomas2004}, to a particular vector field $\bb b:\Sigma\to\R^3$ and
	    suitable surface $\Sigma$ with boundary curve $\partial\Sigma$. We use it in the form
	    \begin{equation}
			\label{eq:stokes}
	     \int_\Sigma\left[(\nabla_\x\times\bb b)\cdot\bs\nu_\Sigma\right](\x)\,\dd S 
			= \int_{\partial\Sigma}[\bb b\cdot\bs{\tau}](\x)\,\dd l.
	    \end{equation}
			Here the normals $\bs\nu$ on the surface and the tangents $\bs\tau$ of the curve are oriented by the right hand rule.
			If you look at the surface from the side to which the normals point, the curve is traversed counter-clockwise.
			
	    Now we take the particular vector field $\bb b = \bs\nu_\Sigma\times\bb a$. Using \eqref{eq:surf_div_tangent}, 
			we obtain for the left hand side of \eqref{eq:stokes}
	    \[
	        \int_\Sigma\left[(\nabla_\x\times(\bs\nu_\Sigma\times\bs a))\cdot\bs\nu_\Sigma\right](\x)\,\dd S
	        = \int_\Sigma\left[\nabla_\Sigma\cdot\bs a^\shortparallel\right](\x)\,\dd S.
			\]
			Using the circular symmetry of the triple product, we have for the right hand side of \eqref{eq:stokes}
			\[
	        \int_{\partial\Sigma}\left[(\bs\nu_\Sigma\times\bb a)\cdot\bs{\tau}\right](\x)\,\dd l 
					= \int_{\partial\Sigma} \left[(\bs{\tau}\times\bs\nu_\Sigma)\cdot\bb a\right](\x)\,\dd l
	        = \int_{\partial\Sigma} \left[\bb a\cdot\bs\nu_{\partial\Sigma}\right](\x)\,\dd l.
	    \]
	    Combining them gives formula \eqref{eq:surface_div}. Equation \eqref{eq:surface_div_2} follows by 
			using \eqref{def:surf_divergence_spacevec}.
	\end{prf}

	Next we obtain an analog of the Reynolds Transport Theorem \ref{thm:reynolds_transport} for moving surfaces.
	 \begin{thm}[Transport Theorem for Moving Surfaces in $\R^3$]
	\label{thm:surface_transport}   
	        Let $\Sigma(\t ), \t  \in [0,T]$ be a family of smooth surfaces in $\Omega\subseteq\mathbb{R}^3$ 
					with the velocity field $\bb w$ of surface points and the surface area element $\dd S = \sqrt{g}\dd u^1\dd u^2$
					for the surface coordinates $\bb u=(u^1, u^2) \in \mathcal{U}$, see Subsection \ref{subsec:surf}.
	        Here $g(\t ,\bb u)$ denotes the determinant of the metric tensor \eqref{def:inverse_det}
					of the surface at time $\t $. 
					
	        Further let $\psi:I\times\Sigma(\t)\to\R$ be a smooth function in the sense that
	        the Lagrangian time derivative $\mathring{\psi}(\t,\x)=\frac D{D\t}\psi(\t,\x)$ exists, see \eqref{def:lagrange_deriv} 
					and \eqref{eq:mueller_deriv}. 
	        Recall that the surface divergence $\nabla_\Sigma\cdot\bb w$ is given by \eqref{def:surf_div}. We set
					$w_\nu= \bb w\cdot\bs\nu_\Sigma$. 
				  We obtain the following transport equations
	        \begin{align}
					\label{eq:surface_transport}
	            \frac{\dd}{\dd \t }\int\limits_{\Sigma(\t )} \psi (\t,\x)\,\dd S
							&=
							\int\limits_{\Sigma(\t )} \left[\mathring{\psi} + \psi\,\nabla_\Sigma\cdot\mathbf{w}\right](\t,\x)\,\dd S
							\\
							\label{eq:surf_trans1}
							&=
							\int\limits_{\Sigma(\t )} \left[\mathring{\psi} + \psi\,\left(\nabla_\Sigma\cdot\mathbf{w}^\shortparallel
							-2\kappa_M w_\nu\right)\right](\t,\x)\,\dd S
							\\
							\label{eq:surf_trans3}
							&=
							\int\limits_{\Sigma(\t )} \left[\frac\delta{\delta\t}\psi 
							+\nabla_\Sigma\cdot \left(\psi\,\mathbf{w}^\shortparallel\right)
							-2\kappa_M w_\nu\right](\t,\x)\,\dd S
							.
					\end{align}
					Now we consider $\psi:I\times\Omega\to\R$ is a smooth function. We we set $\psi_\nu =\bs\nu_\Sigma\cdot\nabla_\x\psi$.
				  Then we can write the transport equations as
	        \begin{align}
							\label{eq:surf_trans2}
							\frac{\dd}{\dd \t }\int\limits_{\Sigma(\t )} \psi (\t,\x)\,\dd S
							&=
							\int\limits_{\Sigma(\t )} \left[\frac{\partial\psi}{\partial \t} + \nabla_\Sigma\cdot(\psi\,\mathbf{w})
							+\psi_\nu\bb w_\nu\right](\t,\x)\,\dd S
							\\
							\label{eq:surf_trans}
							&= 
							\int\limits_{\Sigma(\t )} \left[\frac{\partial\psi}{\partial \t}
							+ \nabla_\Sigma\cdot (\psi\mathbf{w}^\shortparallel) + (\psi_\nu- 2\kappa_M \psi)w_\nu\right] (\t,\x)\,\dd S
							.
	        \end{align}
	    \end{thm}
	\begin{prf}    
	    The proof of the transport theorem for the surface is quite analogous to the proof of Theorem \ref{thm:reynolds_transport}.
	    Important preparations were made in Subsection \ref{subsec:move} by defining the parametrization of the initial
			reference surface $\Sigma_0\subset\Omega$, the flow map $\bs\chi_\Sigma^\t$ 
			as well as the appropriate time derivatives for moving surfaces and clarifying the latter.
	    First we go via the mapping $\bs\Phi_{\Sigma_0}$ from the surface parameters in $\mathcal U$
			to the initial reference configuration $\Sigma_0$ and 
			then via the flow map $\bs\chi_\Sigma^\t=\bs\chi_\Sigma(\t,\cdot)$ to the current surface $\Sigma(\t)$.
	    This mapping is given by $\bs\Phi_\Sigma = \bs\chi_\Sigma^\t\circ\bs\Phi_{\Sigma_0}: \mathcal{U} \to \Sigma(\t )$, 
			see Subsection \ref{subsec:move}. 
			We recall the surface integral given by \eqref{def:surf_int}, 
			the function $\widetilde{\psi}(\t,\bb u)=\psi(\t ,\bs\Phi_\Sigma(\t ,\bb u))$,
			the time derivative of the metric $g$ satisfying \eqref{eq:metricdet_deriv} 
			and \eqref{eq:mueller_deriv} for the Lagrangian derivative. Thus we obtain
	    \begin{align*}
	        \frac{\dd}{\dd \t }\int\limits_{\Sigma(\t )}& \psi(\t,\mathbf{x})\,\dd S
	        = \frac{\dd}{\dd \t }\int\limits_{\mathcal{U}} \psi(\t ,\bs\Phi_\Sigma(\t ,\bb u))\,
					\sqrt{g(\t ,\bb u)}\,\dd u^1\dd u^2\notag\\
					=& \frac{\dd}{\dd \t }\int\limits_{\mathcal{U}} \widetilde{\psi}(\t ,\bb u)\,
					\sqrt{g(\t ,\bb u)}\,\dd u^1\dd u^2
	        = \int\limits_{\mathcal{U}} \left[\frac D{D\t}\psi +
	        \psi\nabla_\Sigma\cdot\mathbf{w}\right](\t ,\bs\Phi_\Sigma(\t ,\bb u))
					\sqrt {g(\t ,\bb u )}\,\dd u^1\dd u^2\\
					=&\int\limits_{\Sigma(\t )}\left(\mathring\psi(\t ,\x) +
	        \psi(\t ,\x)\nabla_\Sigma\cdot\mathbf{w}(\t ,\x)\right)\,\dd S.
	    \end{align*}
			This is the first equation \eqref{eq:surface_transport}.
			The second equation \eqref{eq:surf_trans1} follows using \eqref{def:surf_divergence_spacevec}.
			For the third eqution \eqref{eq:surf_trans3} we use \eqref{eq:lagrange_deriv}.
			The fourth equation \eqref{eq:surf_trans2} is an application of \eqref{eq:lagrange_deriv}, \eqref{eq:thomas_deriv} and
			\eqref{eq:product_rule} to \eqref{eq:surface_transport}. 
			Finally we use \eqref{def:surf_divergence_spacevec} in \eqref{eq:surf_trans2} to obtain \eqref{eq:surf_trans}.   
	\end{prf}

The formula \eqref{eq:surface_transport} can be found in 
Dziuk and Elliott \cite[Lemma 2.2, Appendix A]{Dziuk2007}. But, for the proof they assume that $\psi:I\times\Omega\to\R$.
They called this transport equation a Leibniz formula.
The formula \eqref{eq:surf_trans} is \cite[(2.10)]{Dziuk2007}. They used a level set description of the surface.
The resulting transport equation is the same.
					
The formula \eqref{eq:surf_trans} is also given in Cermelli et al.\ \cite[(3.14)]{Cermelli2005} 
for the case $\psi_\nu =0$ and using \eqref{eq:surface_div} of the 
Divergence Theorem \ref{thm:surface_int}. They used an extension
of $\psi$ to a neighborhood of $\Sigma(\t)$.
See also the literature cited there. The formula \eqref{eq:surf_trans1} 
was given in Petryk and Mroz \cite[(2.34)]{Petryk1986}.
					
In particular note that the time derivative of the metric corresponds to the surface divergence of the velocity field 
$\mathbf{w}$ in \eqref{eq:surface_transport}. This highlights that the surface transport theorem is an analogue to the 
Reynolds Transport Theorem \ref{thm:reynolds_transport}.
The time derivative of the integrated quantity is the integral of the Lagrangian time derivative of a quantity 
plus the quantity times the divergence of the velocity field. The formulas \eqref{eq:surf_trans1}
and \eqref{eq:surf_trans} show how a non-flat geometry comes into play through the mean curvature term.					

Note that Dreyer \cite{Dreyer2003} and M\"uller \cite{Mueller1985} considered only 
the two maps $\bs\Phi_{\Sigma_0}$ and $\bs\Phi_\Sigma$, see Subsection \ref{subsec:move}.
Then $\bs\chi_\Sigma = \bs\Phi_\Sigma (\cdot,\bs\Phi_\Sigma^{-1}(\cdot))$ is only implied.
This leads to slightly more complicated formulas.
They introduced an extra notation for the metric on
$\Sigma_0$ which is actually only $\bb G(0)$ in our notation. For the Transport Theorem \ref{thm:surface_transport}
they refer back to the stationary surface $\Sigma_0$ whereas we refer to its stationary parameter set $\mathcal U$. 
This is mathematically equivalent. It seems that they wanted to stay closer to the proof of the 
Reynolds Transport Theorem \ref{thm:reynolds_transport} for volumes.
					
We find our approach somewhat simpler, as it does not involve the reciprocal of $g(0)$ in the integrals.
Also the introduction of $\bs\chi_\Sigma$ seems to add clarity to the concepts.
Also the formulas \cite[(48)]{Dreyer2003} and 
\cite[(3.9)]{Mueller1985} were given with respect to the parameter variables. 
The time derivative of $\psi$ is then the Lagrangian
time derivative $\mathring{\psi}$ given as \eqref{eq:mueller_deriv}.

Aris \cite[Sections 10.12,10.31]{Aris1989} has a different formulation and proof. He considered a dynamic
function from $\R^2$ to $\R^2$. Both are equipped with Riemanian geometries.
The densities $\psi$ are not functions of $\x\in\R^3$ but of $\bb u(t)\in\R^2$. 
Instead of using \eqref{eq:metricdet_deriv} for the 
differentiated area element he used $\frac{\dd}{\dd\t}\sqrt g=\frac 1{2g}\frac{\dd g}{\dd\t}\, \sqrt g $ instead.
Thereby some insight into the transport is lost. His formulas are completely intrinsic to the surface, which is
mathematically nice, but not so useful for practical applications to the embedded moving surfaces we are interested in. 
In the some sections he did discuss the embedding into $\R^3$, but only in a somewhat sketchy manner. An intrinsic formulation
does not have normal transport of the surface but only the transport within the surface.
		
If we have only intrinsic information on the surface, then there are no normal velocities at all.
So we would also have $w_\nu =0$. Then the equation \eqref{eq:surf_trans} becomes a 2-dimensional analog of
of Reynolds transport theorem \ref{thm:reynolds_transport}. It has a non-Cartesian divergence of the intrinsic
curved geometry of the surface. This case is of interest here only if the surface has a
purely tangential motion.

	\subsection{Separate Differential Balance Equations}
	\label{subsec:sep}

	We consider an open subset $\Omega\subseteq\R^3$ and any family of open control volumes $\mathcal V(\t)$ with 
	$\mathcal V(\t)\subset\Omega$ for $\t\in [0,T]$ with $T>0$. We take	
	the generic balance law \eqref{eq:int_balance_law} on a time dependent control volume $\mathcal V(\t)$. 
	We use formula \eqref{eq:reynolds_thm} of the Reynolds Transport Theorem \ref{thm:reynolds_transport} 
	and \eqref{eq:surf_trans} of the Surface Transport Theorem \ref{thm:surface_transport} 
	to determine the time derivatives. One could take \eqref{eq:surf_trans3} instead with a different time derivative
	and one term less. The latter makes it more general. We assume that the boundary flux $\bb j$ is defined everywhere in the volume
	and $\bb j_\Sigma$ everywhere on the surface. To the volume boundary flux term 
	we apply the usual Gauss-Green theorem and to the surface boundary flux term the surface version 
	\eqref{eq:surface_div} in order to obtain
	\begin{align}
	\label{eq:int_balance}
	0=&
	\tikzmarkin[fill=green,opacity=0.2]{sep:box1}(0,-0.7)(0,0.6)
	\int\limits_{\mathcal{V}(\t )}\left[\frac{\partial}{\partial \t }\psi	  
		+ \nabla_\x\cdot\left[\psi\mathbf{v}+\bb j\right]-\xi\right](\t ,\x)\,\dd V
		\tikzmarkend{sep:box1}\nonumber\\
		&+
		\tikzmarkin[fill=cyan,opacity=0.2]{sep:box2}(0,-0.7)(0,0.6)
		\int\limits_{\Sigma(\t )} \left[\frac{\partial\psi}{\partial \t}
							+ \nabla_\Sigma\cdot (\psi\mathbf{w}^\shortparallel) + (\psi_\nu- 2\kappa_M \psi)w_\nu
		+\nabla_\Sigma\cdot\bb j^\shortparallel_\Sigma 
		-\xi_\Sigma\right](\t,\x)\,\dd S 	
							\tikzmarkend{sep:box2}.
	\end{align}
	This equation is satisfied for any $\t\in [0,T]$ and any control volume $\mathcal V(\t)\subset\Omega$ of any shape and 
	appropriate size. 
	If a control volume does not contain a surface $\Sigma (\t)$ then we take the integral over $\Sigma (\t)$ to be $0$. 
	
	We will make use of the Lebesgue differentiation theorem of real analysis. For this
	we consider for some fixed $\t\in\; ]0,T[$ a family of open balls $B_{\x,r}^n =\{\y\in\R^n\, |\, |\y-\x|<r\}\subset\mathcal V(\t)$ 
	with radius small $r>0$ enough. Their volumes are denoted as $\vol(B_{\x,r}^n)$. 
	Any other shapes with a scaling factor like the radius $r>0$ 
	can be used for such  arguments. 
	
	Take $(\t,\x)\in\; ]0,T[\times\Omega$. Then we can find an open control volume with $\x\in\mathcal V(\t)\subset\Omega$
	and there are infinitely many balls
	$B_{\x,r}^n\subset\mathcal V(\t)$ with $r>0$ sufficiently small.
	Then the Lebesgue differentiation theorem states that for any integrable function $f :\Omega\to\R$ 
	the right hand side of of the inequalities
	\begin{align*}
	\left|f (\x)-\frac 1{\vol(B_{\x,r}^n)}\int_{B_{\x,r}^n}f (\y)\,\dd \y\right|
	&=\left|\frac 1{\vol(B_{\x,r}^n)}\int_{B_{\x,r}^n}[f (\x)-f (\y)]\,\dd \y\right|\\
	&\le	\frac 1{\vol(B_{\x,r}^n)}\int_{B_{\x,r}^n}|f (\x)-f (\y)|\,\dd \y
	\end{align*}
	goes to $0$ for $r\to 0$ for almost any $\x\in\Omega$, see e.g.\ Evans \cite[Appendix E.4]{Evans1998}
	or Rudin \cite[Theorem 7.7]{Rudin1987}. If $f$ is a continuous function this holds for all
	$\x\in\Omega$. This implies that
	\[
	\frac 1{\vol(B_{\x,r}^n)}\int_{B_{\x,r}^n}f (\y)\,\dd \y\;\to\;f(\x)\qquad\text{for}\quad r\to 0.
	\]
	In our application all integrals equal $0$ implying $f(\x)=0$. The variable $\t$ is just a parameter in this argument.
	We assume here that the integrands are continuous. Then the results hold for
	all points $(\t,\x)\in\; ]0,T]\times\Omega$. 
	
	We first consider any $\x_0\in\Sigma(\t)$. We write the surface integral over the open parameter set $\mathcal U$ 
	using \eqref{def:surf_int}. Let $\bb u_0\in\mathcal U$ be such that $\bs\Phi(\bb u_0)=\x_0$. We take a family of balls
	$B_{\bb u_0,r}^2\subset\mathcal U$ for $r>0$ sufficiently small and set $B_{\Sigma,\x_0,r}^2= \bs\Phi[B_{\bb u_0,r}^2]$ 
	to give curved discs on the surface. We define
	the curved {\em pillbox set} over the disc $B_{\Sigma,\x_0,r}^2$ as
	\[
	\mathcal V_{\Sigma,r}=\{\,\x\in\R^3\, |\, \x =\lambda\bs\nu_\Sigma 
	+\y\;\text{with}\; \y\in B_{\Sigma,\x_0,r}^2,\;\lambda\in]-r,r[\, \}
	\subset\Omega.
	\]
	We take the equation \eqref{eq:int_balance} over $\mathcal V_{\Sigma,r}$, then $\Sigma(\t)=B_{\Sigma,\x_0,r}^2$ and
	divide it by $\vol(B_{\bb u_0,r}^2)$. Further,
	we write the volume integral as an average multiplied by the factor $\frac{\vol(V_{\Sigma,r)}}{\vol(B_{\bb u_0,r}^2)}$. 
	The average will 	be bounded and the factor will converge to $0$ for $r\to 0$. This gives the differential equation
	\begin{equation}
	\label{eq:surface_eq}
	 \tikzmarkin[fill=cyan,opacity=0.2]{surf_balance}(0,-0.4)(0,0.6)
	\frac{\partial\psi}{\partial \t}+ \nabla_\Sigma\cdot (\psi\mathbf{w}^\shortparallel) + (\psi_\nu- 2\kappa_M \psi)w_\nu
	+\nabla_\Sigma\cdot\bb j^\shortparallel_\Sigma=\xi_\Sigma
	\tikzmarkend{surf_balance}
	\end{equation}
	satisfied for any $(\t,\x)\in\; ]0,T]\times\Omega$ with $\x\in\Sigma(\t)$ on a surface. 
	Using the integrand of \eqref{eq:surface_transport} it can also be written as
	\[
	\tikzmarkin[fill=cyan,opacity=0.2]{surf_balance1}(0,-0.4)(0,0.6)
	\mathring{\psi}_\Sigma + \psi_\Sigma\,\nabla_\Sigma\cdot\mathbf{w} +\nabla_\Sigma\cdot\bb j^\shortparallel_\Sigma=\xi_\Sigma
	\tikzmarkend{surf_balance1}		
	\]
	in case $\psi_\nu =0$ or $\psi_\nu$ does not exist.
	Now that we know that the integrand of the surface integral vanishes. 
	We can do the analogous argument with averages of the volume integral. 
	This gives us the differential equation for the volume
	\begin{equation}
	\label{eq:volume_eq}
	    \tikzmarkin[fill=green,opacity=0.2]{volume_balance}(0,-0.4)(0,0.6)
	    \frac{\partial}{\partial \t }\psi + \nabla_\x\cdot\left(\psi\mathbf{v} + \bb j \right) 
			= \xi \tikzmarkend{volume_balance}.
	\end{equation}
	for any $(\t,\x)\in\; ]0,T]\times\Omega$. We obtained two uncoupled equations.

	\subsection{A Generalized Transport Theorem}

	As in Subsection \ref{subsec:control}, we want to consider an open control volume 
	$\mathcal V(\t)\subset\Omega\subseteq\R^3$ for $\t\in I=[0,T]$. We assume that it is 
	divided into two open subsets $\mathcal V_1(\t)$ and $\mathcal V_2(\t)$ by an internal open surface $\Sigma(\t)$ with 
	$\mathcal V(\t) =\mathcal V_1(\t)\cup\Sigma(\t)\cup\mathcal V_2(\t)$. This is depicted in Figure \ref{fig:gen_volume}. 
	We also take into account functions $\psi : I\times\Omega \to \mathbb{R}$
	which are continuous in $\mathcal{V}_1(\t)\cup\mathcal{V}_2(\t)$ with limiting values from each side at
	the surface $\Sigma(\t)$. In order to keep the notation compact, we introduce the jump brackets $\dbl \cdot\dbr$.
	For points $\mathbf{x}_\Sigma  \in \Sigma $ we write
	\begin{align*}
	    \psi_i(\t ,\mathbf{x}_\Sigma ) = \lim_{\mathcal{V}_i \ni \mathbf{x} \to \mathbf{x}_\Sigma } \psi(\t ,\mathbf{x}),\, i = 1,2
	    \quad\text{and}\quad
	    \dbl \psi\dbr(\t ,\mathbf{x}_\Sigma ) = \psi_2(\t ,\mathbf{x}_\Sigma ) - \psi_1(\t ,\mathbf{x}_\Sigma ).
	\end{align*}
	The plus sign in the jump bracket is always on the side to which the unit normal vector field $\bs\nu_\Sigma$ on the surface points.
	Usually we will not write out the arguments since either the equations containing jump brackets should 
	hold in any point of the surface or it is clear from the context which point is meant.
	
	Further, we assume that the velocity field $\bb v$ is also discontinuous at $\Sigma$ with bounded limiting values from each side.
	In shock waves, the normal velocity field $\bb v^\perp=\bb v\cdot\bs\nu_\Sigma$ is 
	discontinuous, whereas the 	tangential velocity field $\bb v^\shortparallel$ is continuous, 
	see e.g.\ Anderson \cite[Section 4.3]{Anderson2003} or Toro \cite[Section 3.2]{Toro2009}. 
	The same is true for contact discontinuities, 
	see Anderson \cite[Section 7.1]{Anderson2003} or Toro \cite[p.\ 97]{Toro2009}.
	There may also be shear waves, see Toro \cite[pp.\ 107,111]{Toro2009}, in which the tangential velocity jumps.
	The velocity field has bounded one-sided limits at $\Sigma(\t )$. 
			
	The fact that the velocity field is discontinuous
	at $\Sigma(\t)$ means that the trajectories $\x(\t)=\bs\chi^\t(\x_0)$ are only continuous there and 
	abruptly change direction at an angle. We always assume that the motions $\bs\chi$ are at least continuous
	with respect to $\t$ and $\x$.			
	Let  $\bb v_i$ for $i=1,2$ distinguish the limiting values at $\Sigma(\t)$ of the velocity field on $\mathcal V_i(\t)$.
	The surface $\Sigma(\t )$ moves with the smooth velocity field 
	$\bb w: I\times \Sigma(\t )\to \mathbb{R}^3$
	that may be independent of $\mathbf{v}$. Under these assumptions we want to give a generalization of 
	Theorem \ref{thm:reynolds_transport} and obtain a coupled balance law on the surface.
		
	The proof of Theorem \ref{thm:reynolds_transport} revealed that the total change in time of a physical quantity or state 
	is balanced by the local change of the density function and two additional effects.
	First there is the convective part and second there is a contribution due to the deformation of the material volume.
	We now want to extend it to the case of a singular surface along which the volume functions may be discontinuous.
	\begin{thm}[Generalized Reynolds Transport Theorem in $\R^3$]
	\label{thm:gen_reynolds_transport}
			Let $\mathcal V(\t )\subset\Omega\subseteq\R^3$ for $\t\in [0,T]$ be an open material volume as in Figure \ref{fig:gen_volume}
			with boundary $\mathcal S(\t)=\partial\mathcal V(\t)$. We assume that the boundary $\mathcal S(\t)$
			consists of finitely many smooth pieces. 
			The moving internal two-dimensional surface $\Sigma(\t )$ separates $\mathcal V(\t )$ into two volumes
	    $\mathcal V_1(\t )$ and $\mathcal V_2(\t )$.
			
	    The function $\psi : I\times\Omega \to \mathbb{R}$ and the velocity field 
			$\mathbf{v}: I\times\Omega \to \mathbb{R}^3$ are continuously 
			differentiable on $\mathcal V_1(\t )$ and $\mathcal V_2(\t )$ but discontinuous with bounded one-sided limits
			at the moving surface $\Sigma(\t )$ in the manner introduced above.
	    Then the following generalization of \eqref{eq:reynolds_thm} holds
	    \begin{align}
			\label{eq:gen_reynolds_thm}
	    %
	    \frac{\dd}{\dd \t }\int\limits_{\mathcal V(\t )}\psi(\t ,\mathbf{x}) \,\dd V
			=&
			\int\limits_{\mathcal{V}_1(\t )}\left[\frac{\partial}{\partial \t }\psi
					+ \nabla_\x\cdot(\psi\mathbf{v})\right](\t ,\x) \,\dd V
					+
					\int\limits_{\mathcal{V}_2(\t )}\left[\frac{\partial}{\partial \t }\psi
					+ \nabla_\x\cdot(\psi\mathbf{v})\right](\t ,\x) \,\dd V
					\notag\\
					&-
					\int_{\Sigma(\t)}\Big[\dbl \psi(\mathbf{w}-\bb v)\dbr\cdot\bs{\nu}_\Sigma\Big](\t ,\x)\,\dd S.
			\end{align}
			We can also generalize \eqref{eq:reynolds_thm1} as
			\begin{align}
			\label{eq:gen_reynolds_thm1}
			\frac{\dd}{\dd \t }\int\limits_{\mathcal V(\t )}\psi(\t ,\mathbf{x}) \,\dd V
			=&
			\int\limits_{\mathcal V_1(\t )\cup\mathcal V_2(\t )} \frac{\partial}{\partial \t }\psi(\t ,\mathbf{x})\,\dd V
	    +
			\int\limits_{\mathcal S(\t )} \Big[\psi(\mathbf{v}\cdot\bs\nu)\Big](\t ,\x)\,\dd S
			\notag\\
	    &-
			\int\limits_{\Sigma(\t )}\Big[\dbl \psi\dbr(\mathbf{w}\cdot\bs{\nu}_\Sigma)\Big](\t ,\x)\,\dd S.
			%
	    \end{align}
	\end{thm}
	\begin{prf}
	The integrand on the right hand side of \eqref{eq:reynolds_thm} is singular on $\Sigma(\t)$. The same
	holds for $\Sigma_0$. So in the argument \eqref{eq:reynolds_proof}
	we have to split the integrals over $\mathcal V(\t)$ and $\mathcal V(0)$ into the three parts
	$\widehat{\mathcal V}_1(\t)$, $\mathcal V_D(\t)$, $\mathcal V_2(\t)$ as well as
	$\mathcal V_1(0)$, $\mathcal V_{21}^\t$, $\mathcal V_{22}^\t$. They are defined in Subsection \ref{subsec:control}.
	
	We set $h(\t,\y)=\psi(\t,\bs\chi^\t(\y))\mathcal J^\t(\t,\y)$ for $\y\in\mathcal V(0)$ and obtain
	\begin{align*}
	\frac{\dd}{\dd\t}\left(\int_{\mathcal V_1(\t)}+\int_{\mathcal V_2(\t)}\right)\psi(\t,\x)\,\dd V &=
	\frac{\dd}{\dd\t}\left(\int_{\widehat{\mathcal V}_1(\t)}+\int_{\mathcal V_D(\t)}+\int_{\mathcal V_2(\t)}\right)\psi(\t,\x)\,\dd V\\
	=&\int_{\mathcal V_1(0)}\frac{\partial}{\partial\t}h(\t,\y)\,\dd V_0
	+\frac{\dd}{\dd\t}\left(\int_{\mathcal V_{21}(t)}+\int_{\mathcal V_{22}(t)}\right)h(\t,\y)\,\dd V_0.
	\end{align*}
	The first integral gives the $\widehat{\mathcal V}_1(\t)$ part of the integral over $\mathcal V_1(\t)$ 
	in \eqref{eq:gen_reynolds_thm}.
	
	We assume for the moment that the function $\psi$ and the velocity field $\bb v$ are continuous at $\Sigma$.
	We recall the notations and assumptions from Subsection \ref{subsec:control}.
	Specifically we use the description of the initial volume $\mathcal V_2(0)$
	by \eqref{eq:v2_para} over the open shadow domain $D\subseteq\R^2$.
	The moving surface $\Sigma_0$ in $\mathcal V(0)$ is parametrized using a function $\sigma:I\times D\to\R$.
	Since we are interested in using this theorem locally, this is not a serious restriction. The assumption
	can always be achieved when the surface is small enough. But, the proof can be extended to surfaces that
	do not satisfy the assumption via a piecewise approach by cutting up the surface into suitable parts that do.
	We omit this technical step.
	
	We divide $D$ into two domains, namely $D_1$ where $\sigma\le\alpha$ or $\beta\le\sigma$ and $D_2$ 
	where $\alpha<\sigma<\beta$. The first may be empty. We also introduce the short hand notation
	$\y_{\sigma^i} =(\yin_1,\yin_2,\sigma^i(\t,\yin_1,\yin_2))$ with $i=1,2$ 
	in order to indicate the one-sided limiting values of $h$ on $\Sigma$.
	For the moment we assume that they are identical, but we keep track of them.
	We can now use the Leibniz integral rule with respect to $\yin_3$ to obtain
	\begin{align*}
	\frac{\dd}{\dd\t}\left(\int_{\mathcal V_{21}(t)}+\int_{\mathcal V_{22}(t)}\right)&h(\t,\y)\,\dd V_0
	=\int_{D_1}\left(\int_\alpha^\beta \frac{\dd}{\dd\t}h(\t,\y)\,\dd\yin_3\right)\,\dd A\\
	&+\int_{D_2}\frac{\dd}{\dd\t}\left(\int_{\alpha}^{\sigma}h(\t,\y)\,\dd\yin_3\right)\,\dd A
	+\int_{D_2}\frac{\dd}{\dd\t}\left(\int_{\sigma}^{\beta}h(\t,\y)\,\dd\yin_3\right)\,\dd A\\
	=&\int_D\int_{\alpha}^{\beta}\frac{\partial}{\partial\t}h(\t,\y)\,\dd\yin_3\,\dd A
	+\int_{D_2}\left[h(\y_{\sigma^1})\sigma^1_\t-h(\y_{\sigma^2})\sigma^2_\t\right]\,\dd A\\
	=&\left(\int_{\mathcal V_{21}(t)}+\int_{\mathcal V_{22}(t)}\right)\frac{\partial}{\partial\t}h(\t,\y)\,\dd V_0
	+\int_{D_2}\left[h(\y_{\sigma^1})\sigma^1_\t-h(\y_{\sigma^2})\sigma^2_\t\right]\,\dd A.
	\end{align*}
	The volume integrals give the remaining contributions to volume integrals in \eqref{eq:gen_reynolds_thm}.
	For the surface integral we have using
	\eqref{eq:sigma_t}, \eqref{def:surf_int} and \eqref{eq:sigma_0_int}
	\[
	\int_{D_2}h(\y_{\sigma^1})\sigma^1_\t\,\dd A 
	=\int_{D_2}\psi(\t,\bs\chi^\t(\y_{\sigma^1}))\mathcal J^\t(\t,\y_{\sigma^i})\,
	u_\nu(\t,\y_{\sigma^1})|\bs\tau_1^0\times\bs\tau_2^0|(\t,\y_{\sigma^1})\,\dd A.
	\]
Now we apply \eqref{eq:u_def} and \eqref{eq:int_surf_trans} to obtain
	\begin{align*}
	\int_{D_2}h(\y_{\sigma^1})\sigma^1_\t\,\dd A&\\
	=&\int_{D_2}\psi(\t,\bs\chi^\t(\y_{\sigma^1}))\mathcal J^\t(\t,\y_{\sigma^i})
	\left(\left(\bb D_\x\,\bs\chi^{-\t}(\t,\x)\right)(\bb w(\t,\x)-\bb v^i(\t,\x))\right)\\
	&\qquad\qquad\qquad\qquad\qquad\qquad\qquad\qquad
	\cdot\bs\nu_{\Sigma_0}(\t,\y_{\sigma^1})|\bs\tau_1^0\times\bs\tau_2^0|(\t,\y_{\sigma^1})\,\dd A\\
	=&\int_{D_2}\left[\psi(\bb w-\bb v^1)\cdot\bs\nu_\Sigma
	|\bs\tau_1\times\bs\tau_2|\right](\t,\bs\chi^\t(\y_{\sigma^1}))\,\dd A\\
	=&\int_\Sigma \left[\psi(\bb w-\bb v^1)\cdot\bs\nu_\Sigma\right](\t,\x)\,\dd S.
	\end{align*}
	This leads to \eqref{eq:gen_reynolds_thm} since we obtain
	\[
	\int_{D_2}\left[h(\y_{\sigma^1})\sigma^1_\t-h(\y_{\sigma^2})\sigma^2_\t\right]\,\dd A=
	\int_\Sigma \left[\,\dbl\psi(\bb w-\bb v)\dbr\cdot\bs\nu_\Sigma\right](\t,\x)\,\dd S.
	\]
	
	Since we assumed that the function $\psi$ and the velocity field $\bb v$ are continuous at $\Sigma$ the latter integral vanishes. 
	But the formula \eqref{eq:gen_reynolds_thm} remains true in the discontinuous case 
	since we never used derivatives of $\psi$ and $\bb v$ on $\Sigma$ in its derivation, but only their
	one-sided limits.
	
	We can apply the Gauss-Green Theorem to the integrals over $\mathcal V_1(\t)$ and $\mathcal V_2(\t)$. 
	On $\Sigma$ this gives an integrand $\dbl\psi\bb v\dbr$ that cancels out leaving the velocity field
	$\bb w$ to give \eqref{eq:gen_reynolds_thm1}.
	\end{prf}

	This extension of the Reynolds Transport Theorem \ref{thm:reynolds_transport} is found in Thomas \cite{Thomas1949}, 
	Truesdell and Toupin \cite[Section 192]{Truesdell1960} and M\"uller \cite[Section 3.1]{Mueller1985} with rather heuristic proofs.
	There it is stated in the form of formula \eqref{eq:gen_reynolds_thm1} which differs from \eqref{eq:gen_reynolds_thm} 
	only by an application of the Gauss-Green theorem.

	\subsection{Derivation of Surface Equation with Jump Conditions}
	\label{subsec:jump}

	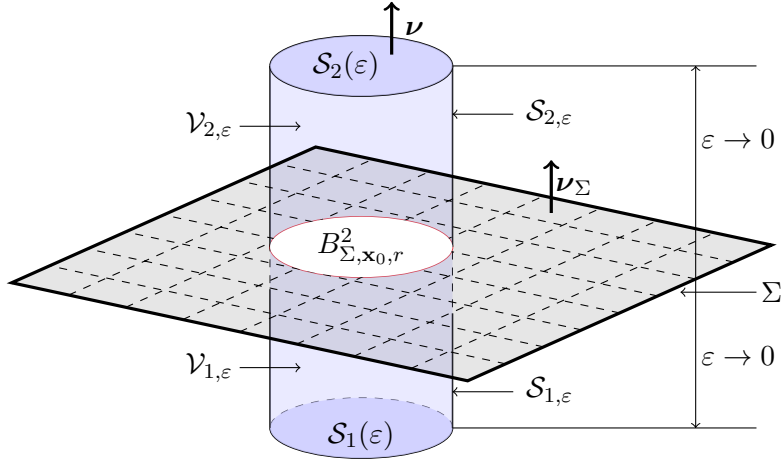
\begin{figure}[!ht]
	    \centering
	    %
	    \begin{tikzpicture}
	        \begin{scope}[scale = 0.8,rotate=90,xshift=2cm,yshift=-1cm]
	            %
	            \draw[dashed] (0,0) arc (-90:90:0.5 and 1.5);
	            \draw[semithick] (0,0) arc (270:90:0.5 and 1.5);
	            \fill[blue!20,fill opacity=0.8] (0,0) arc (-90:90:0.5 and 1.5) arc (90:270:0.5 and 1.5) -- cycle;
	            %
	            \draw[semithick,red] (3,0) arc (-90:90:0.5 and 1.5);
	            \draw[semithick,red] (3,0) arc (270:90:0.5 and 1.5);
	            %
	            \draw[semithick] (0,0) -- (1.9,0);
	            \draw[semithick] (0,3) -- (2.3,3);
	            \draw[dashed] (2,0) -- (3,0);
	            \draw[dashed] (2.4,3) -- (3,3);
	            \draw[semithick] (3,0) -- (6,0);
	            \draw[semithick] (3,3) -- (6,3);
	            \fill[blue!20, fill opacity=0.4] (0,0) arc (270:90:0.5 and 1.5) to (6,3) arc (90:270:0.5 and 1.5) to (0,0) -- cycle;
	            \draw[->] (5,3.5) -- (5,2.5);
	            \draw (5,4) node {$\mathcal{V}_{2,\varepsilon}$};
	            \draw[->] (1,3.5) -- (1,2.5);
	            \draw (1,4) node {$\mathcal{V}_{1,\varepsilon}$};
							\draw [->] (5.2,-1) -- (5.2,0);
					    \draw (5.2,-1.6) node {$\mathcal S_{2,\varepsilon}$};
							\draw [->] (0.6,-1) -- (0.6,0);
					    \draw (0.6,-1.6) node {$\mathcal S_{1,\varepsilon}$};
	            %
	            \draw[semithick] (6,0) arc (-90:90:0.5 and 1.5);
	            \draw[semithick] (6,0) arc (270:90:0.5 and 1.5);
	            \fill[blue!20,fill opacity=0.8] (6,0) arc (-90:90:0.5 and 1.5) arc (90:270:0.5 and 1.5) -- cycle;
	            %
	            \draw (6,0) -- (6,-5);
	            \draw (0,0) -- (0,-5);
	            \draw [<->] (6,-4) -- (0,-4);
							\draw (4.8,-4.7) node {$\varepsilon\to 0$};
	            \draw (1.2,-4.7) node {$\varepsilon\to 0$};
	       \end{scope}
	       \begin{scope}[xshift=0.5cm,yshift=3cm,every node/.append style={yslant=0.5,xslant=-1},yslant=0.45,xslant=-1.5]
	            %
	            \fill[black,fill opacity=0.1] (-1,-1) rectangle (3,3);
	            \draw[step=5mm,black,dashed] (-1,-1) grid (3,3);
	            \draw[black, very thick] (-1,-1) rectangle (3,3);
	            %
	        \end{scope}
					\draw [very thick,->] (0,6.55) -- (0,7.25);
	        \draw (0.3,6.9) node {$\bs{\nu}$};
	        \draw [very thick,->] (2.1,4.45) -- (2.1,5.15);
	        \draw (2.4,4.8) node {$\bs{\nu}_\Sigma$};
					\draw [->] (4.8,3.4) -- (3.8,3.4);
	        \draw (5,3.4) node {$\Sigma $};
					\draw (-0.4,1.5) node {$\mathcal S_1(\varepsilon)$};
					\draw (-0.6,6.4) node {$\mathcal S_2(\varepsilon)$};
					\begin{scope}[color=white]
					\fill (-0.4,4) circle [x radius=1.2, y radius=0.4];
					\end{scope}
					\draw (-0.4,4) node {$B^2_{\Sigma,\x_0,r}$};
	    \end{tikzpicture}
	    \caption{Sketch of the pillbox argument.}
	    \label{fig:pillbox}
	\end{figure}

	To obtain the balance equation for the singular points on a surface there is a so called {\em pillbox argument}. It was
	used by Truesdell and Toupin \cite[Section 193]{Truesdell1960}, M\"uller \cite[Section 3.1]{Mueller1985} 
	and for the two dimensional case by Gurtin \cite[Section 3.2]{Gurtin1993}. The result generalizes the derivation
	of jump conditions on surfaces of discontinuity that goes back to the seminal paper of Riemann \cite{Riemann1990}.
	Zempl\'en \cite{Zemplen1905} contributed a more formal and general derivation of jump conditions for quasi one dimensional flows in 
	cylindrical ducts. This was then extended by Kotchine \cite{Kotchine1926} to more general surfaces. He used pillbox
	type arguments and took jumps in the tangential velocity into account.
	
	We consider an open subset $\Omega\subseteq\R^3$ and any surface of singular points 
	$\Sigma(\t)\subset\Omega$ for $\t\in [0,T]$.
	The idea is to choose a family of cylindrical type sets $\mathcal V_{\Sigma,\varepsilon}\subset\Omega$ 
	as neighborhoods around the surface as control volume.
	We define analogous to Subsection \ref{subsec:sep} a new curved {\em pillbox set}, for any $\varepsilon >0$
	small enough, over the curved disc $B_{\Sigma,\x_0,r}^2$ for any $\x_0\in\Sigma(\t)$ and some suitable $r>0$ as
	\[
	\mathcal V_{\Sigma,\varepsilon}=\{\,\x\in\R^3\, |\, \x 
	=\lambda\bs\nu_\Sigma +\y\;\text{with}\; \y\in B_{\Sigma,\x_0,r}^2,\;\lambda\in]-\varepsilon,\varepsilon[\, \}
	\subset\Omega
	\]
	with the bottom surface $\mathcal S_1(\varepsilon)$ and top surfaces $\mathcal S_2(\varepsilon)$  given as
	\[
	\mathcal S_i(\varepsilon)=\{\,\x\in\R^3\, |\, \x 
	=(-1)^i\varepsilon\bs\nu_\Sigma +\y\;\text{with}\; \y\in B_{\Sigma,\x_0,r}^2\, \}
	\]
	as well as side surfaces 
	\[
	\mathcal S_{i,\varepsilon}=\{\,\x\in\R^3\, |\, \x 
	=(-1)^i\lambda\bs\nu_\Sigma +\y\;\text{with}\; \y\in \partial B_{\Sigma,\x_0,r}^2,\;\lambda\in]0,\varepsilon[\, \}.
	\]
	Here we do not want to shrink the curved disc on the surface, as we did in Subsection \ref{subsec:sep}, but only the normal direction.
	We divide $\mathcal V_{\Sigma,\varepsilon}$ into $\mathcal V_{1,\varepsilon}$, defined as the subset with
	$\lambda\in ]-\varepsilon,0[$ and analogously $\mathcal V_{2,\varepsilon}$ 
	with $\lambda\in ]0,\varepsilon[$. The unit normal vector field $\bs\nu_\Sigma$ on the
	surface $\Sigma$ points into $\mathcal V_{2,\varepsilon}$.
	A sketch for a flat surface is given in Figure \ref{fig:pillbox}.
	
	We consider the general balance law \eqref{eq:int_balance_law} on such an arbitrary material control volume 
	$\mathcal V_{\Sigma,\varepsilon}$. Additionally, we assume that the flux $\bb j$ is discontinuous on $\Sigma$
	with bounded one-sided limits. For the time derivative of the volume terms we use \eqref{eq:gen_reynolds_thm1} 
	and for the surface term \eqref{eq:surf_trans} to obtain
	\begin{align*}
	    \tikzmarkin[fill=green,opacity=0.2]{pill:box1}(0,-0.7)(0,0.6)
	    \int\limits_{\mathcal V_{1,\varepsilon}\cup\mathcal V_{2,\varepsilon}}&\frac{\partial}{\partial \t }\psi\,\dd V
			\tikzmarkend{pill:box1}
			+	\int\limits_{\mathcal S_{1,\varepsilon}\cup\mathcal S_{2,\varepsilon}} \psi(\mathbf{v}\cdot\bs\nu)\,\dd S
			+ 
			\tikzmarkin[fill=cyan,opacity=0.2]{pill:box2}(0,-0.7)(0,0.6)
			\int\limits_{\mathcal S_1(\varepsilon)\cup\mathcal S_2(\varepsilon)} \psi(\mathbf{v}\cdot\bs\nu)\,\dd S
			 \tikzmarkend{pill:box2}
	    -
			\tikzmarkin[fill=cyan,opacity=0.2]{pill:box2a}(0,-0.9)(0,0.6)
			\int\limits_{B_{\Sigma,\x_0,r}^2}\dbl \psi\dbr(\mathbf{w}\cdot\bs{\nu}_\Sigma)\,\dd S
	    \tikzmarkend{pill:box2a}\\
	    &\qquad+ \tikzmarkin[fill=cyan,opacity=0.2]{pill:box3}(0,-0.8)(0,0.6)
	    \int\limits_{B_{\Sigma,\x_0,r}^2} \left [\frac\partial{\partial\t}\psi_\Sigma
			+ \nabla_\Sigma\cdot (\psi_\Sigma\mathbf{w}^\shortparallel) + (\psi_\nu- 2\kappa_M \psi_\Sigma)w_\nu\right]\,\dd S
	    \tikzmarkend{pill:box3}
			\\
	    =& -\tikzmarkin[fill=cyan,opacity=0.2]{pill:box4}(0,-0.7)(0,0.6)
	    \int\limits_{\mathcal S_{1,\varepsilon}\cup\mathcal S_{2,\varepsilon}\cup
			\mathcal S_1(\varepsilon)\cup\mathcal S_2(\varepsilon)} \bb j\cdot\bs\nu\,\dd S
	    \tikzmarkend{pill:box4}
	    %
	    -\tikzmarkin[fill=cyan,opacity=0.2]{eq41:box4}(0,-0.8)(0,0.6)
	    \int\limits_{\partial B_{\Sigma,\x_0,r}^2} \bb j_\Sigma\cdot\bs\nu_{\partial\Sigma}\,\dd l
	    \tikzmarkend{eq41:box4}
	    + \tikzmarkin[fill=green,opacity=0.2]{eq41:box5}(0,-0.7)(0,0.6)
	    \int\limits_{\mathcal V_{1,\varepsilon}\cup\mathcal V_{2,\varepsilon}}\xi\,\dd V
	    \tikzmarkend{eq41:box5}
	    +\tikzmarkin[fill=cyan,opacity=0.2]{eq41:box6}(0,-0.8)(0,0.6)
	    \int\limits_{B_{\Sigma,\x_0,r}^2} \xi_\Sigma\,\dd S.
	    \tikzmarkend{eq41:box6}
	\end{align*}
	Note that here $\psi_\nu=\nabla_\x\psi_\Sigma\cdot\bs\nu_\Sigma$ and $w_\nu =\bb w\cdot\bs\nu_\Sigma$.
	In case $\psi_\nu =0$ or $\psi_\nu$ does not exist, we would replace the use of \eqref{eq:surf_trans} by \eqref{eq:surf_trans3}.
	Now as $\varepsilon$ goes to zero, the volume integrals over $\mathcal V_{1,\varepsilon}$ and $\mathcal V_{2,\varepsilon}$
	as well as the surface integrals over $\mathcal S_{1,\varepsilon}$ and $\mathcal S_{2,\varepsilon}$ go to $0$. 
	These are the green and white terms.
	The surface integrals over $\mathcal S_1(\varepsilon)$ and $\mathcal S_2(\varepsilon)$ become integrals over $B_{\Sigma,\x_0,r}^2$
	where for $\mathcal S_1(\varepsilon)$ we have $\bs\nu \to -\bs\nu_\Sigma$ 
	and $\bs\nu \to\bs\nu_\Sigma$ for $\mathcal S_2(\varepsilon)$.
	The integrals connected to $B_{\Sigma,\x_0,r}^2$ remain unchanged. But, to the integral over $\partial B_{\Sigma,\x_0,r}^2$
	we apply \eqref{eq:surface_div}. This gives
	\begin{align*}
	    0=&\int\limits_{B_{\Sigma,\x_0,r}^2} \left[\left(\dbl\psi\mathbf{v}\dbr
	    - \dbl \psi\dbr\mathbf{w}\right)\cdot\bs{\nu}_\Sigma
	    + \frac{\partial}{\partial \t }\psi_\Sigma
			+ \nabla_\Sigma\cdot (\psi_\Sigma\mathbf{w}^\shortparallel) + (\psi_\nu- 2\kappa_M \psi_\Sigma)w_\nu\right](\t,\x)\,\dd S\\
	    & +\int\limits_{B_{\Sigma,\x_0,r}^2} \left[\dbl\bb j\dbr\cdot\bs\nu_\Sigma\
	    +\nabla_\Sigma\cdot \bb j_\Sigma^\shortparallel -\xi_\Sigma\right](\t,\x)\,\dd S.
	\end{align*}

	Now we use the argument from Subsection \ref{subsec:sep}, since $\x_0$ and $r>0$ can be arbitrarily chosen,
	that the equality holds point-wise at $(\t,\x)$ for the integrands. Thus we obtain the point-wise balance law on the surface
	\begin{equation}
	\label{eq:gen_loc_jump_cond}
	\tikzmarkin[fill=cyan,opacity=0.2]{sing_balance}(0,-0.4)(0,0.6)
	\frac{\partial}{\partial \t }\psi_\Sigma
			+ \nabla_\Sigma\cdot (\psi_\Sigma\mathbf{w}^\shortparallel) + (\psi_\nu- 2\kappa_M \psi_\Sigma)w_\nu
			+\nabla_\Sigma\cdot \bb j_\Sigma^\shortparallel -\xi_\Sigma
			=\dbl \psi\dbr w_\nu -\dbl\psi\mathbf{v}+\bb j\dbr\cdot\bs\nu_\Sigma	
	\tikzmarkend{sing_balance}.				
	\end{equation}
	On the right hand side is the usual jump condition for conservation laws. The left hand side
	consists of the extra surface terms that we included.
	
	In case we have $\psi_\Sigma$ only defined on the surfaces $\Sigma$ we can use \eqref{eq:surface_transport}
	in the balance law \eqref{eq:int_balance_law} and follow the same steps. We can replace \eqref{eq:gen_loc_jump_cond} by
	the rather concise form
	\begin{equation}
	\tikzmarkin[fill=cyan,opacity=0.2]{sing_balance1}(0,-0.4)(0,0.6)
	\label{eq:gen_loc_jump_cond1}
	\mathring{\psi}_\Sigma + \psi_\Sigma\,\nabla_\Sigma\cdot\mathbf{w}
			+\nabla_\Sigma\cdot \bb j_\Sigma^\shortparallel -\xi_\Sigma
			=\dbl \psi\dbr w_\nu -\dbl\psi\mathbf{v}+\bb j\dbr\cdot\bs\nu_\Sigma
				\tikzmarkend{sing_balance1}.
	\end{equation}
	We may substitute $\nabla_\Sigma\cdot\bb w=\nabla_\Sigma\cdot\bb w^\shortparallel -2\kappa_Mw_\nu$ 
	from \eqref{def:surf_divergence_spacevec}. Further we use 
	\eqref{def:surf_divergence_spcevec}, \eqref{eq:decomp} and \eqref{def:surface_metric} to have
	$\nabla_\Sigma\cdot \mathbf{w}^\shortparallel=\nabla_\alpha w_\tau^\alpha$ and 
	$\nabla_\Sigma\cdot \bb j_\Sigma^\shortparallel=\nabla_\alpha  j_\tau^\alpha$. This gives these equations 
	in the form
	\begin{equation}
	\label{eq:gen_loc_jump_cond2}
	\tikzmarkin[fill=cyan,opacity=0.2]{sing_balance2}(0,-0.4)(0,0.6)
	\mathring{\psi}_\Sigma + \psi_\Sigma(\nabla_\alpha w_\tau^\alpha- 2\kappa_Mw_\nu) 
			+\nabla_\alpha j_\tau^\alpha -\xi_\Sigma
			=\dbl \psi\dbr w_\nu -\dbl\psi\mathbf{v}+\bb j\dbr\cdot\bs\nu_\Sigma	
	\tikzmarkend{sing_balance2}.
	\end{equation}
	This is the formula in M\"uller \cite[(3.16)]{Mueller1985} in the variables $(\t,\bb u)$.
	There the source term $\xi_\Sigma$ is split into a production and a supply term, see also
	Mavrovouniotis and Brenner \cite[(5.1)]{Mavrovouniotis1993}. 
	Dreyer \cite[(51)]{Dreyer2003} used $\frac\partial{\partial\t}\widetilde{\psi}_\Sigma(\t,\bb u)$ instead
	of the identical $\mathring{\psi}_\Sigma(\t,\bb u)$, see \eqref{eq:mueller_deriv}.

	\section{Moving Curves in Two Dimensions, Points in One Dimension}

	So far we derived the equations in $\R^3$. But the lower dimensional cases 
	may also be derived from these equations. In particular we focus on (\ref{eq:gen_loc_jump_cond}).
	We briefly want to {\em sketch} how this may be done in $\R^2$. Details of the theory of moving 
	piecewise smooth curves in the plane
	together with physical properties are presented in Gurtin \cite{Gurtin1993}.
	
	We assume a symmetry such that the surface density is constant in one surface direction 
	and that the normal vector fields have no component in this direction.
	For example we may choose our coordinate system such that the density is constant in 
	$\xin_3$ direction and the normal vector fields have the $\xin_3$ component zero.
	Then the problem reduces to a two dimensional problem, where the moving singular surface $\Sigma(\t)$ reduces 
	to a moving curve $\c(\t)$
	in the $\xin_1$-$\xin_2$ plane. The remaining surface parameter may be taken to be the {\em arc-length} $\s$ of the curve.
	Then the tangent vector is a unit vector. Further the mean curvature, as the arithmetic mean of the principle curvatures, 
	reduces to the curvature $\kappa$ of the curves.
	Finally the divergence of the tangential component of the surface velocity reduces to the 
	derivative with respect to the arc-length $\s$. The Lagrangian derivative $\mathring{\psi}$
	is an advective tangential derivative with respect to arc-length. 
	Thus we obtain as a modification of \eqref{eq:gen_loc_jump_cond} a formula, for the case in which 
	$\psi_c$ is also defined in the plane, as
	\begin{align}
	\label{eq:gen_loc_jump_cond_2d}
	    \left[\mathring{\psi}_\c  + \psi_\c \left(\frac{\partial}{\partial \s}w_\tau 
			- \kappa w_\nu\right)+\psi_{\c,\nu}w_\nu- \xi_\c + \frac{\partial}{\partial \s}\bb j_\c^\shortparallel\right](\t ,\s)
	    = \left[w_\nu\dbl \psi\dbr - \dbl \psi\mathbf{v} + \bb j \dbr\bs{\nu}\right] (\t ,\s) 
	\end{align}
	If $\psi_\c$ is defined only on the moving curves we can
	use the formula without the normal derivative of $\psi_\c$, i.e.\ $\psi_{\c,\nu}=0$.
	
	Now for the one dimensional case we assume a planar interface, which stays planar in time, with a spatially constant surface density.
	We have a moving point of singularity $\p(\t)\in\R$.
	Since the metric is constant in time, the derivative $\partial_\t {g}$ vanishes and there is no curvature.
	The velocity is the scalar function $w(\t,\p(\t))= \frac\dd{\dd\t}\p(\t)$. Thus we obtain
	\begin{align}
	    \frac{\dd}{\dd \t }\psi_\p (\t,\p(\t))  - \xi_p(\t,\p(\t))= \left[w\dbl \psi\dbr 
			- \dbl \psi v +  j\dbr\right](\t,\p(\t))  .
			\label{eq:gen_loc_jump_cond_1d}
	\end{align}
	This type of jump condition is also known as {\em generalized Rankine-Hugoniot jump conditions}, cf.\ Yang \cite[(3.7)]{Yang1999}.
	\appendix
	\section{Time Dependent Determinants}
	\label{app:det}

We want to calculate the time derivative of the determinant of the Jacobian matrix or
deformation gradient of the motion $\bs{\chi }^{t}$
from Subsection \ref{subsec:control}. It is denoted by ${\mathcal J}^{\t }$ , i.e.\
\begin{equation}
\label{def:J}
{\mathcal J}^{\t } =	\, \det \bb D_\x  \, \bs \chi  ^\t  =
\; \det \left( \frac{\partial
\bs{\chi }^{\t }_{j}}{\partial x _{k}}\right)_{1\leq j,\, k \leq n } .
\end{equation}
For this let $\bb A = (a_{jk})_{1\leq j,\, k \leq n}\in\R^{n\times n}$
a matrix or tensor of order 2 with $\t$ dependent coefficients. For any pair of indices
$j,k$ we define the matrix $\bb A_{jk}$ as the matrix obtained from $\bb A$
by replacing row $j$ and column $k$ zero entries except for $a_{jk}=1$, i.e.\
$$
  \bb A_{jk} =
  \left(
    \begin{array}{cccccc}
      a_{11} & \cdots & 0 & \cdots & a_{1n } \\
      \vdots &	      &   &	   & \vdots \\
      0      & \cdots & 1 & \cdots & 0      \\
      \vdots &	      &   &	   & \vdots \\
      a_{n 1} & \cdots & 0 & \cdots & a_{n n }
    \end{array}
  \right) .
$$
With the {\em co-factors} $\det\bb A_{jk}$ of $\bb A$ the inverse matrix $\bb A^{-1}$is given as, 
see Meyer\cite[Section 6.2]{Meyer2001} or Strang 
\cite[Section 4.4]{Strang1988}\footnote{These references do not use the matrix $\bb A_{jk}$ but the
$(n-1)\times(n-1)$ where the $j$-th row and $k$-th column are deleted. The determinant then involves an extra factor
of $(-1)^{j+k}$ that we can avoid.},
$$
\bb A^{-1} = \left( a^{-1}_{jk} \right) _{1\leq j,\, k \leq n }= \left( a^{jk} \right) _{1\leq j,\, k \leq n }
= \left( \frac{\det \bb A_{kj}}
{\det \bb A}\right)_{1 \leq j,\, k \leq n } ,
$$
i.e.\ we have
\begin{equation}
  \label{eq:1.9}
  \det \bb A_{jk}= a^{-1}_{kj} \, \det \bb A= a^{kj} \, \det \bb A
\end{equation}
and $ a_{jl} a^{lk} = \delta_j^k$. We use the summation convention that indices that appear twice are
summed from 1 to $n$ and $\delta_{jk}=\delta_j^k = 0$ for $j \neq k$ and $=1$ for $j=k$
are the {\em Kronecker delta} tensors. We introduce the $n$th order {\em permutation tensor} $\epsilon_{i_1\cdots i_n}$
using indexed indices. It
is equal to $1$ when $i_1\cdots i_n$ is an even permutation of the numbers $123\cdots n$, equal to $-1$ for
an odd permutation and otherwise equal to $0$. The determinant of $\bb A$ is
$$
\det \bb A = \epsilon_{i_1\cdots i_n}  a_{1 i_1} \cdots  a_{n i_n }
$$
and and that of $\bb A_{jk} $
$$
\det \bb A_{jk} = \epsilon_{i_1\cdots i_n} a_{1i_1} \cdots  a_{j-1\, i_{j-1}}\cdot \delta_{i_jk}\cdot
a_{j+1\, i_{j+1}}\cdots a_{ni_n}.
$$
This implies that $\frac{\partial \det \bb A}{\partial a_{jk}}= \det \bb A_{jk}$.
Now with (\ref{eq:1.9}) we have
\begin{equation}
\label{eq:mat_det_deriv}
\frac{d \det \bb A}{dt}  = \frac{\partial \det \bb A}{\partial a_{jk}}\,
\frac{\partial a_{jk}}{dt} = \det \bb A_{jk}\,\frac{da_{jk}}{dt} = a^{kj} \,\frac{da_{jk}}{dt}\,\det \bb A .
\end{equation}
This we apply to ${\mathcal J}^{\t }$ and note that in this case
$a^{-1}_{kj}= a^{kj}=\frac{\partial x_{k}}{\partial \bs\chi ^{\t }_{j}}$. We obtain
\[
\frac{d{\mathcal J}^{\t }}{d\t }  = \frac{\partial x_{k}}{\partial \bs\chi ^{\t }_{j}}\,
\frac{d}{d\t } \left(\frac{\partial \bs\chi ^{\t }_{j}}{\partial x_{k}}\right)\,{\mathcal J}^{\t }
= \frac{\partial x_{k}}{\partial \bs\chi ^{\t }_{j}}\,\frac{\partial v^j}{\partial x_{k}}\, {\mathcal J}^{\t } 
= \frac{\partial v^j}{\partial \bs\chi ^{\t }_{j}}\,{\mathcal J}^{\t }
= \frac{\partial v^j}{\partial x_{j}}\,{\mathcal J}^{\t }
\]
or
\begin{equation}
\label{eq:det_deriv}
  \frac{d{\mathcal J}^{\t }}{d\t } = \left(  \nabla_\x \cdot  \, \bb {v} \right)\,{\mathcal J}^{\t } .
\end{equation}
The two dimensional analogue is given in \eqref{eq:metricdet_deriv}.

\section{Cross Product Transformation}

Let $\bb A\in\R^{3\times 3}$ be an invertible matrix
and $\bs\tau_1,\bs\tau_2\in\R^3$ two vectors. We denote by $\bb C =\left(\det\bb A_{jk}\right)_{1\le j,k\le 3}$ the
matrix of co-factors. Then we have the following formula for the cross product.
\begin{equation}
\label{eq:prod_trans}
\bb A\bs\tau_1\times\bb A\bs\tau_2 =\bb C (\bs\tau_1\times\bs\tau_2)=\det\bb A(\bb A^{-1})^T(\bs\tau_1\times\bs\tau_2).
\end{equation}
From \eqref{eq:1.9} we have $\bb A^{-1} =\frac 1{\det\bb A}\bb C^T$. This gives the second equality.
The first equality is a slightly tedious calculation for each of the three components. 
Starting from the left hand side some terms cancel to give the right hand side.

    %
    \bibliographystyle{abbrv}
    \bibliography{literatur_jc}
\end{document}